\documentclass[11pt]{article}
\usepackage{amsfonts,mathrsfs,amssymb,amsthm,mathptm}
\usepackage{amsmath,amscd}
\usepackage{mathptm,pslatex}
\usepackage{color}
\usepackage[dvips]{graphicx}
\usepackage[all]{xy}
\usepackage{graphicx}
\usepackage{fancyhdr}

\begin{document}
\newtheorem{Def}{Definition}[section]
\newtheorem{Bsp}[Def]{Example}
\newtheorem{Prop}[Def]{Proposition}
\newtheorem{Theo}[Def]{Theorem}
\newtheorem{Lem}[Def]{Lemma}
\newtheorem{Koro}[Def]{Corollary}
\theoremstyle{definition}
\newtheorem{Rem}[Def]{Remark}

\newcommand{\add}{{\rm add}}
\newcommand{\con}{{\rm con}}
\newcommand{\gd}{{\rm gl.dim}}
\newcommand{\sd}{{\rm st.dim}}
\newcommand{\sr}{{\rm sr}}
\newcommand{\dm}{{\rm dom.dim}}
\newcommand{\cdm}{{\rm codomdim}}
\newcommand{\tdim}{{\rm dim}}
\newcommand{\E}{{\rm E}}
\newcommand{\Mor}{{\rm Morph}}
\newcommand{\End}{{\rm End}}
\newcommand{\ind}{{\rm ind}}
\newcommand{\rsd}{{\rm res.dim}}
\newcommand{\rd} {{\rm rd}}
\newcommand{\ol}{\overline}
\newcommand{\overpr}{$\hfill\square$}
\newcommand{\rad}{{\rm rad}}
\newcommand{\soc}{{\rm soc}}
\renewcommand{\top}{{\rm top}}
\newcommand{\pd}{{\rm pdim}}
\newcommand{\id}{{\rm idim}}
\newcommand{\fld}{{\rm fdim}}
\newcommand{\Fac}{{\rm Fac}}
\newcommand{\Gen}{{\rm Gen}}
\newcommand{\fd} {{\rm fin.dim}}
\newcommand{\Fd} {{\rm Fin.dim}}
\newcommand{\Pf}[1]{{\mathscr P}^{<\infty}(#1)}
\newcommand{\DTr}{{\rm DTr}}
\newcommand{\cpx}[1]{#1^{\bullet}}
\newcommand{\D}[1]{{\mathscr D}(#1)}
\newcommand{\Dz}[1]{{\mathscr D}^+(#1)}
\newcommand{\Df}[1]{{\mathscr D}^-(#1)}
\newcommand{\Db}[1]{{\mathscr D}^b(#1)}
\newcommand{\C}[1]{{\mathscr C}(#1)}
\newcommand{\Cz}[1]{{\mathscr C}^+(#1)}
\newcommand{\Cf}[1]{{\mathscr C}^-(#1)}
\newcommand{\Cb}[1]{{\mathscr C}^b(#1)}
\newcommand{\Dc}[1]{{\mathscr D}^c(#1)}
\newcommand{\K}[1]{{\mathscr K}(#1)}
\newcommand{\Kz}[1]{{\mathscr K}^+(#1)}
\newcommand{\Kf}[1]{{\mathscr  K}^-(#1)}
\newcommand{\Kb}[1]{{\mathscr K}^b(#1)}
\newcommand{\DF}[1]{{\mathscr D}_F(#1)}

\newcommand{\Kac}[1]{{\mathscr K}_{\rm ac}(#1)}
\newcommand{\Keac}[1]{{\mathscr K}_{\mbox{\rm e-ac}}(#1)}

\newcommand{\modcat}{\ensuremath{\mbox{{\rm -mod}}}}
\newcommand{\Modcat}{\ensuremath{\mbox{{\rm -Mod}}}}
\newcommand{\Spec}{{\rm Spec}}

\newcommand{\stmc}[1]{#1\mbox{{\rm -{\underline{mod}}}}}
\newcommand{\Stmc}[1]{#1\mbox{{\rm -{\underline{Mod}}}}}
\newcommand{\prj}[1]{#1\mbox{{\rm -proj}}}
\newcommand{\inj}[1]{#1\mbox{{\rm -inj}}}
\newcommand{\Prj}[1]{#1\mbox{{\rm -Proj}}}
\newcommand{\Inj}[1]{#1\mbox{{\rm -Inj}}}
\newcommand{\PI}[1]{#1\mbox{{\rm -Prinj}}}
\newcommand{\GP}[1]{#1\mbox{{\rm -GProj}}}
\newcommand{\GI}[1]{#1\mbox{{\rm -GInj}}}
\newcommand{\gp}[1]{#1\mbox{{\rm -Gproj}}}
\newcommand{\gi}[1]{#1\mbox{{\rm -Ginj}}}

\newcommand{\opp}{^{\rm op}}
\newcommand{\otimesL}{\otimes^{\rm\mathbb L}}
\newcommand{\rHom}{{\rm\mathbb R}{\rm Hom}\,}
\newcommand{\pdim}{\pd}
\newcommand{\Hom}{{\rm Hom}}
\newcommand{\Coker}{{\rm Coker}}
\newcommand{ \Ker  }{{\rm Ker}}
\newcommand{ \Cone }{{\rm Con}}
\newcommand{ \Img  }{{\rm Im}}
\newcommand{\Ext}{{\rm Ext}}
\newcommand{\StHom}{{\rm \underline{Hom}}}
\newcommand{\StEnd}{{\rm \underline{End}}}
\newcommand{\KK}{I\!\!K}
\newcommand{\gm}{{\rm _{\Gamma_M}}}
\newcommand{\gmr}{{\rm _{\Gamma_M^R}}}

\def\vez{\varepsilon}\def\bz{\bigoplus}  \def\sz {\oplus}
\def\epa{\xrightarrow} \def\inja{\hookrightarrow}

\newcommand{\lra}{\longrightarrow}
\newcommand{\llra}{\longleftarrow}
\newcommand{\lraf}[1]{\stackrel{#1}{\lra}}
\newcommand{\llaf}[1]{\stackrel{#1}{\llra}}
\newcommand{\ra}{\rightarrow}
\newcommand{\dk}{{\rm dim_{_{k}}}}

\newcommand{\holim}{{\rm Holim}}
\newcommand{\hocolim}{{\rm Hocolim}}
\newcommand{\colim}{{\rm colim\, }}
\newcommand{\limt}{{\rm lim\, }}
\newcommand{\Add}{{\rm Add }}
\newcommand{\Prod}{{\rm Prod }}
\newcommand{\Tor}{{\rm Tor}}
\newcommand{\Cogen}{{\rm Cogen}}
\newcommand{\Tria}{{\rm Tria}}
\newcommand{\Loc}{{\rm Loc}}
\newcommand{\Coloc}{{\rm Coloc}}
\newcommand{\tria}{{\rm tria}}
\newcommand{\Con}{{\rm Con}}
\newcommand{\Thick}{{\rm Thick}}
\newcommand{\thick}{{\rm thick}}
\newcommand{\Sum}{{\rm Sum}}

{\Large \bf
\begin{center}
Derived and stable equivalences of centralizer matrix algebras
\end{center}}

\medskip
\centerline{\textbf{Xiaogang Li} and \textbf{Changchang Xi}$^*$ }

\renewcommand{\thefootnote}{\alph{footnote}}
\setcounter{footnote}{-1} \footnote{ $^*$ Corresponding author.
Email: xicc@cnu.edu.cn; Fax: 0086 10 68903637.}
\renewcommand{\thefootnote}{\alph{footnote}}
\setcounter{footnote}{-1}
\renewcommand{\thefootnote}{\alph{footnote}}
\setcounter{footnote}{-1}

\begin{abstract}
The centralizer of a matrix in a full matrix algebra is called a principal centralizer matrix algebra. Characterizations are presented for principal centralizer matrix algebras to be Morita equivalent, almost $\nu$-stable derived equivalent, derived equivalent, and stably equivalent of Morita type, respectively, in terms of new equivalence relations on square matrices. These equivalence relations on matrices are introduced in a natural way by their elementary divisors. Thus the categorical equivalences are reduced to questions in linear algebra. Consequently,  principal centralizer matrix algebras of permutation matrices are Morita equivalent if and only if they are derived equivalent. Moreover, two representation-finite, principal centralizer matrix algebras over a perfect field are stably equivalent of Morita type if and only if they are stably equivalent. Further, derived equivalences between the principal centralizer matrix algebras of permutation matrices induce the ones of their $p$-regular parts and $p$-singular parts of the given permutations.
\end{abstract}

{\footnotesize\tableofcontents\label{contents}}

\section{Introduction\label{Introduction}}
Let $R$ be a field and $n$ a natural number. We denote by $[n]$ the set of integers $\{1,\cdots,n\}$ and by $M_n(R)$ the full $n\times n$ matrix algebra over $R$ with the identity matrix $I_n$.
For a nonempty set $X$ of $M_n(R)$, the centralizer matrix algebra $S_n(X,R)$ of $X$ in $M_n(R)$ is defined by $$S_n(X,R):=\{a\in M_n(R)\mid ax=xa,\; \forall \; x\in X\}.$$

In case of $X=\{c\}$,  we write $S_n(c,R)$ for $S_n(X,R)$. Clearly, $S_n(X,R)=\cap_{c\in X}S_n(c,R)$. The $R$-algebra $S_n(c,R)$ is termed as a \emph{principal centralizer matrix algebra}.

\medskip
Clearly, if $X$ consists of invertible matrices, then the centralizer matrix algebras $S_n(X,R)$ are a special class of invariant algebras, which can be dated back to the classical invariant theory (see \cite{HW}). If $X$ consists of nilpotent matrices over an algebraically closed field $R$, then all nilpotent matrices in $S_n(X,R)$ form a variety which is of significant interest in semisimple Lie algebras (see \cite{DP,AP}). Further, the centralizer matrix algebras $S_n(X,R)$ also appear in the representation theory of finite groups. For instance, if $R$ is a field of characteristic $p>0$ and $G$ is a finite group with a Sylow $p$-subgroup $P$, then $G$ acts on the set $G\backslash P$ of left cosets of $P$ in $G$ by left multiplication, and the set $R[G\backslash P]$ of all $R$-linear combination of the left cosets of $P$ in $G$ is a permutation module over the group algebra $R[G]$ of $G$. Moreover, $\End_{R[G]}(R[G\backslash P])$ is isomorphic to a centralizer matrix algebra $S_r(G_0,R)$, where $r=|G:P|$  and $G_0$ is a subgroup of permutation matrices in $M_r(R)$. In the literature, $\End_{R[G]}(R[G\backslash P])$  is also called  a modular Hecke algebra (see \cite{Al2, ca}).
Note that $G_0$ is isomorphic to $G/O_p(G)$, where $O_p(G)$ stands for the largest normal $p$-subgroup of $G$. Alperin suggested to study the endomorphism ring ${\End}_{R[G]}(R[G\backslash P])$ for attacking the famous Alperin's weight conjecture (see \cite{Al}).

Centralizer matrix algebras have been studied in various aspects of mathematics, such as invariant subspaces or orbits in \cite{bf}, and maximal doubly stochastic matrices in \cite{cdfk}. Recently, in a series of papers \cite{xz1,xz2,xz3}, a lot of new structural and homological properties of $S_n(c,R)$ is revealed. For instance, $S_n(c,R)$ is always a cellular $R$-algebra if $R$ is an algebraically closed field,  and the famous Auslander-Reiten (or Auslander-Alperin) conjecture on stable equivalences holds true for $S_n(c,R)$ over an arbitrary field $R$. The conjecture states that stably equivalent algebras should have the same number of non-projective non-isomorphic simple modules. Further, $S_n(c,R)$ is always a Gorenstein algebra and captures the Auslander algebra of the truncated polynomial algebra $R[x]/(x^n)$ for all $n\in \mathbb{N}$, which has played an important role in the classification of  parabolic subgroups of classical groups with a finite number of orbits on the unipotent radical (see \cite{hr}).

In this note, we study relations between principal centralizer matrix algebras. Particularly, we are interested in when they are Morita, derived and stably equivalent of Marita type. These equivalences are of great importance in the representation theory of algebras and groups \cite{rouq}. Our question here can be formulated precisely as follows.

\medskip
{\bf Question}: Let $R$ be a field, $c\in M_n(R)$ and $d\in M_m(R)$.  What are the necessary and sufficient conditions for $S_n(c,R)$ and $S_m(d,R)$ to be Morita, derived or stably equivalent (of Morita type)?

\smallskip
The answer to this question is closely related to the minimal polynomials of $c$ and $d$, and our characterizations of Morita, derived and stable equivalences are given surprisingly in a very elementary way, namely in terms of elementary divisors of matrices $c$ and $d$.

\smallskip
To state our main result precisely, we first introduce several equivalence relations on square matrices.

Let $R[x]$ be the polynomial algebra over a field $R$ in one variable $x$.
Given polynomials $f(x)$ and  $g(x)$ of positive degree, we define $f(x)\le g(x)$ if $f(x)$ divides $g(x)$, that is, $g(x)=f(x)h(x)$ with $h(x)\in R[x]$.

For $c\in M_n(R)$, let

$\mathcal{E}_c\subset R[x]$ denote the set of distinct elementary divisors of $c$. Here we understand that the elements of a set are pairwise different, and the ones of a multiset are allowed to be duplicate. 

$\mathcal{M}_c:=\{f(x)\in \mathcal{E}_c\mid f(x) \mbox{ is maximal with respect to the order } \le \},$ the set of maximal divisors of $c$.

For $f(x)\in \mathcal{M}_c$, we define

$P_c(f(x)):=\{i\ge 1 \mid  \exists \mbox{ irreducible polynomial } p(x) \mbox{ such that } p(x) \le f(x), p(x)^i\in \mathcal{E}_c \},$ the set of power indices of $f(x)$ in $\mathcal{E}_c$.

$\mathcal{R}_c:=\{f(x)\in \mathcal{M}_c\mid f(x) \mbox{ is a reducible polynomial} \},$ the set of reducible maximal divisors of $c.$

The polynomials in $\mathcal{R}_c$ are exactly those $f(x)\in \mathcal{M}_c$ with $P_c(f(x))\neq \{1\}$.

Let $\mathbb{Z}_{>0}$ be the set of all positive integers and $s\in \mathbb{Z}_{>0}$. For a subset  $T:=\{ m_1, m_2, \cdots, m_s\}$ of $\mathbb{Z}_{>0}$ with $m_1>m_2>\cdots >m_s$, we define a multiset $\mathcal{H}_T:= \{\{m_1-m_2,\cdots,m_{s-1}-m_s,m_s\}\}$ and a set $\mathcal{J}_T=:\{m_1,m_1-m_2,\cdots,m_1-m_s\}$. If $T=\{m_1\}$, then $\mathcal{H}_T=\mathcal{J}_T=T$.

Now we introduce a few new equivalence relations on square matrices.
\begin{Def}
Two matrices $c\in M_n(R)$ and $d\in M_m(R)$ are said to be

$(1)$ $M$-equivalent, written $c\stackrel{M}\sim d$, if there is a bijection $\pi$ between $\mathcal{M}_c$ and $\mathcal{M}_d$ such that $R[x]/(f(x))\simeq R[x]/((f(x))\pi)$ as algebras and ${P_c(f(x))}= {P_d((f(x))\pi)}$ for all $f(x)\in \mathcal{M}_c$, where $(f(x))\pi$ denotes the image of $f(x)$ under the map $\pi$.

$(2)$  $D$-equivalent, written $c\stackrel{D}\sim d$, if there is a bijection $\pi$ between $\mathcal{M}_c$ and $\mathcal{M}_d$ such that $R[x]/(f(x))\simeq R[x]/((f(x))\pi)$ as algebras and $\mathcal{H}_{P_c(f(x))}= \mathcal{H}_{P_d((f(x))\pi)}$ for all $f(x)\in \mathcal{M}_c$.

$(3)$ $AD$-equivalent, written $c\stackrel{AD}\sim d$, if there is a bijection $\pi$ between $\mathcal{M}_c$ and $\mathcal{M}_d$ such that $R[x]/(f(x))\simeq R[x]/((f(x))\pi)$ as algebras and $P_c(f(x))= {P_d((f(x))\pi)}$ or $P_c(f(x))=\mathcal{J}_{P_d((f(x))\pi)}$ for all $f(x)\in \mathcal{M}_c$.

$(4)$ $SM$-equivalent, written $c\stackrel{SM}\sim d$, if there is a bijection $\pi$ between $\mathcal{R}_c$ and $\mathcal{R}_d$ such that  $R[x]/(f(x))\simeq R[x]/((f(x))\pi)$ as algebras and ${P_c(f(x))}={P_d((f(x))\pi)}$ or $P_c(f(x))=\mathcal{J}_{P_d((f(x))\pi)}$ for all $f(x)\in \mathcal{R}_c$.
\end{Def}

Next, we characterize Morita, derived and stable equivalences between principal centralizer matrix algebras in terms of these equivalence relations on matrices. Our main result reads as follows.

\begin{Theo}\label{main1}
Let $R$ be a field, $c\in M_n(R)$ and $d\in  M_m(R)$. Then

$(1)$ $S_n(c,R)$ and $S_m(d,R)$ are Morita equivalent if and only if $c\stackrel{M}\sim d$.

$(2)$ $S_n(c,R)$ and $S_m(d,R)$ are derived equivalent if and only if $c\stackrel{D}\sim d$.

$(3)$ $S_n(c,R)$ and $S_m(d,R)$ are almost $\nu$-stable derived equivalent if and only if $c\stackrel{AD}\sim d$.

$(4)$ Assume that either $R$ is perfect or both $c$ and $d$ are invertible matrices of finite order. Then $S_n(c,R)$ and $S_m(d,R)$ are stably equivalent of Morita type if and only if $c\stackrel{SM}\sim d$.
\end{Theo}
Thus the existence of a Morita equivalence, an almost $\nu$-stable derived equivalence and a derived equivalence between principal centralizer matrix algebras can be read off from the elementary divisors of given matrices directly, and therefore is a problem in linear algebra.

As a consequence, we have the following corollary.

\begin{Koro}\label{derp}Let $R$ be a field, $c\in M_n(R)$ and $d\in  M_m(R)$.

$(1)$ If $c$ and $d$ are permutation matrices, then $S_n(c,R)$ and $S_m(d,R)$ are Morita equivalent if and only if they are derived equivalent.

$(2)$ Suppose that $R$ is perfect and that $S_n(c,R)$ and $S_m(d,R)$ are representation-finite. Then $S_n(c,R)$ and $S_m(d,R)$ are stably equivalent of Morita type if and only if they are stably equivalent.

$(3)$ If $S_n(c,R)$ and $S_m(d,R)$ are derived equivalent or stably equivalent, then they have the same dominant dimension.
\end{Koro}

Moreover, for the centralizer matrix algebras of permutation matrices, a derived equivalence between them gives rise to derived equivalences of smaller centralizer matrix algebras corresponding to $p$-regular and $p$-singular parts of permutations. For details, we refer to Proposition \ref{regular-singular}.

The paper is organized as follows. In Section \ref{sect2} we fix notation, recall basic definitions and terminologies, and prove a few preliminary lemmas needed in the later proofs. In Section \ref{Pf} we prove the main result and its corollary. In Section \ref{sect5} we present examples to show that even for principal centralizer matrix algebras over a field, the notions of  Morita equivalences, almost $\nu$-stable derived equivalences and derived equivalences are different, though they may coincide in some cases. Finally, we propose some open questions for further investigation of centralizer matrix algebras. For example, how to characterize generally the stable equivalences of centralizer matrix algebras?

\section{Preliminaries}\label{sect2}

In this section we recall some basic definitions and terminologies, and prepare a few lemmas for our proofs.

\subsection{Derived equivalences of algebras}
In this paper, $R$ is a field unless stated otherwise. By an algebra we mean a finite-dimensional unitary associative algebra over $R$.

Let $A$ be an algebra. By $\rad(A)$ and $LL(A)$ we denote the Jacobson radical and the Loewy length of $A$, respectively. Let $A^{\opp}$ and $A^e$ stand for the opposite algebra and the enveloping algebra $A\otimes_R A\opp$ of $A$, respectively. By a module we always mean a left module unless stated otherwise. We write $A\modcat$ for the category of all finitely generated left $A$-modules, $A\modcat_{\mathscr{P}}$ for the full subcategory of $A\modcat$ consisting of modules without any nonzero projective summands, and $A\prj$ (respectively, $A\inj$) for the full subcategory of $A\modcat$ consisting of projective (respectively, injective) $A$-modules.

For an $A$-module $M$, we denote by ${\add}(M)$ the full subcategory of $A\modcat$ consisting of all modules isomorphic to direct summands of direct sums of finitely many copies of $M$. Let $\ell(M)$ stand for the composition length of $M$, $\mathcal{B}(M)$ for the basic module of $M$, and $M_{\mathscr{P}}$ for the submodule of $M$ without any nonzero projective summand such that $M/M\!_{\mathscr{P}}$ is projective.

The stable module category of $A$ is denoted by $A\stmc$ which has the same objects as $A$-mod does, but the morphism set $\underline{\Hom}_A(X,Y)$ of objects $X$ and $Y$ is the quotient of $\Hom_A(X,Y)$ modulo $\mathcal{P}(X,Y)$, the set of all homomorphisms that factorize through projective $A$-modues. If $M$ is a non-projective indecomposable $A$-module, then $\mathcal{P}(M,M)\subseteq \rad(\End_A(M))$.

When we say the number of a specific class of $A$-modules such as projective, injective, simple and indecomposable $A$-modules, we always refer to the number of isomorphism classes of them.

For homomorphisms $f:X\to Y$ and $g: Y\to Z$ in $A\modcat$, we write $fg$ for their composition. This implies that the image of an element $x\in X$ under $f$ is denoted by $(x)f$. Thus $\Hom_A(X,Y)$ is naturally a left $\End_A(X)$- and right $\End_A(Y)$-bimodule.  The composition of functors between categories is written from right to left, that is, for two functors
$F:\mathcal{C}\ra \mathcal{D}$ and $G:\mathcal{D}\ra \Sigma$, we write $G\circ F$, or simply $GF$, for the composition of $F$ with $G$. The image of an object $X\in \mathcal{C}$ under $F$ is written as $F(X).$

Let $D: A\modcat\to A^{\opp}\modcat$ be the usual duality of the algebra $A$. The Nakayama functor $\nu_A:=D\Hom_A(-,A)\simeq D(A)\otimes_A-$ from $A\modcat$ to itself restricts to an equivalence: $\prj{A}\stackrel{\sim}{\to} A\inj$. An $A$-module $M$ is said to be $\nu$-\emph{stably projective} if $\nu^i_A{M}$ is projective for all $i\ge 0$. For example, if
$e^2=e\in A$ satisfies $\add(\nu_A{Ae})=\add(Ae)$, then $Ae$ is $\nu$-stably projective. In this case,  $e\in A$ is said to be \emph{$\nu$-stable}. Let $A$-stp denote the full subcategory of $A\modcat$ consisting of all $\nu$-stably projective $A$-modules.
The \emph{Frobenius part} of $A$ is defined (up to Morita equivalence) to be the algebra $eAe$ where $e$ is an idempotent such that $\add(Ae)=A$-stp (see \cite{hx2}).

Let $\D{A}$ (respectively, $\Db{A}$) be the unbounded (respectively, bounded) derived category of $A\modcat$. They are $R$-linear, triangulated categories. Algebras $A$ and $B$ are said to be \emph{derived equivalent} if their derived categories $\Db{A}$ and $\Db{B}$ are equivalent as $R$-linear triangulated categories. An $R$-linear triangle equivalence $F:\Db{A}\rightarrow\Db{B}$ is called a \emph{derived equivalence} between $A$ and $B$.

A special class of derived equivalences, called almost $\nu$-stable derived equivalences, was introduced in \cite{hx1} to establish relations between derived equivalences and stable equivalences of Morita type (see \cite{hx1} for more details). One of the interesting properties of almost $\nu$-stable derived equivalences is that such a derived equivalence between finite-dimensional algebras always induces a stable equivalence of Morita type (see \cite[Theorem 1.1]{hx1}), and thus preserves global and dominant dimensions of algebras. Recall that finite-dimensional algebras $A$ and $B$ are \emph{stably equivalent of Morita type} \cite{MB} if there exist left-right projective (projective as left and right module) bimodules $_AM_B$ and $_BN_A$ such that $M\otimes_B N\simeq A\oplus P$ as $A^e$-modules for some projective $A^e$-module $P$ and $N\otimes_A M\simeq B\oplus Q$ as $B^e$-modules for some projective $B^e$-module $Q$.
Clearly, the exact functor $N\otimes_A-: A\modcat{}\ra B\modcat$ induces a stable equivalence $N\otimes_A-:A\stmc{}\ra B\stmc$.

It is known that $R$ and any separable $R$-algebra are stably equivalent of Morita type. Recall that an $R$-algebra $A$ is separable over $R$ if $A$ is a projective module over $A^e$.

Algebras $A$ and $B$ are said to be {\it Morita equivalent}
if $A\modcat$ and $B\modcat$ are equivalent as $R$-linear categories. Clearly, Morita equivalences  are almost $\nu$-stable derived equivalences.

To get derived equivalences, the following corollary of ${\rm \cite[Theorem 1.1]{hx2}}$ provides a convenient way. For further information on constructions of derived equivalences of algebras, we refer to \cite{x3}.

Let $\mathcal{C}$ be an additive category and $\mathcal{D}$ a full subcategory of $\mathcal{C}$. For $Y\in \mathcal{C},$ a morphism $f:M\ra Y$ with $M\in \mathcal{D}$ is called a \emph{right $\mathcal{D}$-approximation} of $Y$ if each morphism $D\to Y$ with $D\in \mathcal{D}$ factorizes through $f$. Dually, one defines a \emph{left $\mathcal{D}$-approximation} of an object $X$ in $\mathcal{C}.$

A sequence $X\stackrel{g}\ra M\stackrel{f}\ra Y$ in $\mathcal{C}$ with $M\in \mathcal{D}$ is called a \emph{$\mathcal{D}$-split sequence} \cite{hx2} if $g$ is both a kernel of $f$ and a left $\mathcal{D}$-approximation of $X$, and if $f$ is both a cokernel of $g$ and a right $\mathcal{D}$-approximation of $Y$.

\begin{Lem}{\rm\cite{hx2}}
\label{split-thm}
Let $A$ be an algebra, and let $\mathcal{C}$ be a full subcategory of $A\modcat$ and $M$  an object in $\mathcal{C}$. Suppose
$X\ra M'\ra Y$
is an ${\add}(M)$-split sequence in $\mathcal{C}$. Then the endomorphism algebras ${\End}_{\mathcal{C}}(M\oplus X)$ and
 ${\End}_{\mathcal{C}}(M\oplus Y)$ are derived equivalent via a tilting module.
\end{Lem}

The next simple observation characterizes Morita equivalences.
\begin{Lem}\label{M}
Let $A$ be an algebra and $M,N\in A\modcat$. Then  ${\End}_A(M)$ and ${\End}_A(N)$ are Morita equivalent if and only if ${\add}(M)$ and ${\add}(N)$ are equivalent as $R$-linear categories.
\end{Lem}

\begin{Lem}\label{alm}
{\rm \cite[Section 3, Remark]{hx1}}
Let $A$ be a self-injective algebra and $X\in A\modcat$. Then the endomorphism algebras ${\End}_A(A\oplus X)$ and ${\End}_A(A\oplus \Omega_A(X))$ are almost $\nu$-stable derived equivalent, where $\Omega_A(X)$ stands for the syzygy of $X$.
\end{Lem}

\begin{Lem}\label{almst}
{\rm \cite[Theorem 4.4]{CM}}
Let $A$ and $B$ be symmetric algebras, and let $F$ be an almost $\nu$-stable derived equivalence between two gendo-symmetric algebras ${\rm End}_A(A\oplus M)$ and ${\rm End}_B(B\oplus N)$, where $M$ and $N$ are basic non-zero modules without projective summands. Then $A$ and $B$ are (almost $\nu$-stable) derived equivalent. Furthermore, $F$ induces a stable equivalence $\overline{F}:A\stmc$ $\ra B\stmc$ with $\overline{F}(M) = N$.
\end{Lem}

\begin{Lem}\label{iso}
Let $A$ and $B$ be commutative self-injective algebras, $_AM$ and $_BN$ be faithful modules over $A$ and $B$, respectively. If the endomorphism algebras ${\End}_A(M)$ and ${\End}_B(N)$ are derived equivalent, then $A\simeq Z(\End_A(M))\simeq Z\big({\End}_B(N))\simeq B$, where $Z(C)$ denotes the center of an algebra $C$.
\end{Lem}

{\it Proof.} For an algebra $C$ and a faithful $C$-module $X$, one always has an embedding $Z(C)\hookrightarrow Z(\End_C(X))$. Thus $A\hookrightarrow Z(\End_A(M))$ since $A$ is commutative. Note that a faithful module over a self-injective algebra is clearly a generator-cogenerator. This implies that $M_{\End_A(M)}$ is a right faithful module and the bimodule $_AM_{\End_A(M)}$ has the double centralizer property. Thus there is an embedding $Z(\End_A(M))\hookrightarrow \End_{\End_A(M)\opp}(M)\simeq A.$ Hence $A\simeq Z(\End_A(M))$.
Now, assume that ${\End}_A(M)$ and ${\End}_B(N)$ are derived equivalent. Then $Z\big({\End}_A(M)\big)\simeq Z\big({\End}_B(N)\big)$ by \cite[Proposition 9.2]{Rickard1}, and therefore $A\simeq Z(\End_A(M))\simeq Z\big({\End}_B(N))\simeq B$. $\square$

\subsection{Modules over quotients of polynomial algebras}
In this section we remind basic facts on modules over the polynomial algebra $R[x]$.
Throughout this section, $R$ is a field unless stated otherwise.

Let $f(x)$ be an irreducible polynomial in $R[x]$ and $A:=R[x]/(f(x)^n)$ for $n\in \mathbb{N}$. Then $A$ is a self-injective algebra by \cite[Corollary 4.37]{rotman}. Further, $A$ is a local, commutative, symmetric, Nakayama algebra (see, for instance \cite[Example, p.127]{ARS}). Thus $A$ has $n$ indecomposable modules $M_{f(x)}(i):=R[x]/(f(x)^i)$ for $i\in [n].$ For simplicity, we often write $M(i)$ for $M_{f(x)}(i)$ and understand $M(0)=0$. Clearly, $\Hom_A(M(i),A)\simeq \Hom_R(M(i),R)\simeq M(i)$ as $A$-modules for all $i\in [n].$

Let $B:=R[x]/(f(x)^m)$ for $m<n$. Then there is a canonical surjective homomorphism $\pi:A\ra B$ of $R$-algebras, and each $B$-module can be viewed as an $A$-module via $\pi$. Up to isomorphism, indecomposable $A$-modules coming from $B$-modules are exactly those $M(i)$ with $i\in [m].$ For $M, N\in B\modcat$, ${\Hom}_A(M,N)={\Hom}_B(M, N)$.

For an irreducible polynomial $g(x)\in R[x]$ and a positive integers $m$, if $A\simeq R[x]/(g(x)^m)$ as $R$-algebras, then $n=LL(R[x]/(f(x)^n))=LL(R[x]/(g(x)^m))=m$, and for $t\in [n]$, the indecomposable $R[x]/(g(x)^n)$-module $R[x]/(g(x)^t)$ is isomorphic to the $A$-module $R[x]/(f(x)^t)$.

\begin{Lem}\label{lift-exact}
Let $a,b,c,d \in \{0,1,\cdots,n\}$ such that $b < a< c$, $b< d< c$ and $a+d=b+c$. If $_AX\in A\modcat$ has no indecomposable direct summands $N$ with $b<\ell(N)<c$ and $_AY:={}_AX\oplus M(b)\oplus M(c)$,  then there is an ${\add}(_AY)$-split sequence $0\ra M(a)\ra M(b)\oplus M(c)\ra M(d)\ra 0$.
\end{Lem}

{\it Proof.}
Let $g:M(b)\to M(d)$ and $h: M(c)\to M(d)$ be the canonical injective and surjective homomorphisms, respectively. We define
$v:=\left(\begin{smallmatrix} g\\	h\\		\end{smallmatrix} \right).$
Then $v:M (b)\oplus M(c)\to M(d)$ is a surjective homomorphism. Similarly, let $p: M(a)\to M(b)$ and $q: M(a)\to M(c)$ be the canonical surjective and injective homomorphisms, respectively, and $u:=(-p,q)$. Then $u: M(a)\to M(b)\oplus M(c)$ is an injective homomorphism. By the definition of $M(i)$, we have $uv=0.$ It follows from $a+d=b+c$ that the sequence
$$ (\star)\quad 0\lra M(a)\stackrel{u}\lra M(b)\oplus M(c)\lraf{v} M(d)\lra 0$$
of $A$-modules is exact. We shall show that $u$ and $v$ are left and right ${\add}(_AY)$-approximations of $M(a)$ and $M(d)$, respectively. In fact, we need only to show that $v$ is a right ${\add}(_AY)$-approximation of $M(d)$ because applying the dual functor ${\Hom}_R(-,R)$ shows that $u$ is a left ${\add}(_AY)$-approximation of $M(a)$. To show that $v$ is a right ${\add}(_AY)$-approximation of $M(d)$, it suffices to prove that any homomorphism $h:Z\to M(d)$ with $Z$ an indecomposable summand of $Y$ factorizes through $v$. By assumption, either $\ell(Z)\leq b$ or $\ell(Z)\geq c.$ Suppose $\ell(Z)\leq b.$ Then $\ell((Z)h)\leq b=\ell((M(b))g)$, and therefore $(Z)h\subset (M(b))g.$ Let $s: Z\to M(b)\oplus M(c)$ be the homomorphism defined by $(z)s:=(((z)h)g^{-1},0)$ for $z\in Z.$ Clearly, $h=sv$. Suppose $\ell(Z)\geq c.$ Then max$\{a,b,c, d\}\le\ell(Z)$. Let $B:=R[x]/(f(x)^{\ell(Z)})$. Then $B$ is the quotient of $A$ by the ideal $(f(x))^{n-\ell(Z)}$, $Z\simeq M(\ell(Z))=R[x]/(f(x)^{\ell(Z)})=B$ as $A$-modules, and the exact sequence $(\star)$ can be viewed as the one of $B$-modules. So the exactness of ${\Hom}_{B}(Z,-)$ implies that $h$ factorizes through $v$ in $B\modcat$. Since ${\Hom}_A(M, N)={\Hom}_B(M,N)$ for $M, N\in B\modcat$, $h$ factorizes through $v$ in $A\modcat$. $\square$

\begin{Lem}\label{de}
Let $n=\sum^s_{i=1} \ell_i$ with $\ell_i\in \mathbb{Z}_{>0}$. For $\sigma\in \Sigma_s$ and $j\in [s]$, define $M_j:=M(\sum^j_{i=1} \ell _i)$ and  $M^\sigma_j:= M(\sum^j_{i=1} \ell_{(i)\sigma})$. Then the endomorphism algebras ${\End}_A(\bigoplus^s_{j=1} M_j)$ and ${\End}_A(\bigoplus^s_{j=1} M^\sigma_j)$ are derived equivalent.
\end{Lem}

{\it Proof.} The symmetric group $\Sigma_s$ is generated by all transpositions $(t,t+1),t\in [s-1]$. In particular, $\sigma\in \Sigma_s$ can be written as a product of those transpositions, say $\sigma=\prod^k_{i=1} (t_i,t_i+1)$ for $t_i\in [s-1].$ Set $\sigma_{k+1}:=1$ and $\sigma_r:=\prod^k_{i=r} (t_i,t_i+1)$ for all $r\in [k].$ Then $(t_r,t_r+1)\sigma_r=\sigma_{r+1}$.

It suffices to show that there is a derived equivalence between ${\End}_A\big(\bigoplus^s_{j=1} M^{\sigma_r}_j\big)$ and ${\End}_A\big(\bigoplus^s_{j=1} M^{\sigma_{r+1}}_j\big)$ for all $r\in [k].$ For any $\tau\in \Sigma_s$, we define $\sum^{t_r-1}_{i=1} \ell_{(i)\tau}=0$ if $t_r=1$. For $r\in [k]$, let $a_r=\ell_{(t_r+1)\sigma_{r+1}}+\sum^{t_r-1}_{i=1} \ell_{(i)\sigma_{r+1}},b_r=\sum^{t_r-1}_{i=1} \ell_{(i)\sigma_{r+1}},c_r=\sum^{t_r+1}_{i=1} \ell_{(i)\sigma_{r+1}},d_r=\sum^{t_r}_{i=1} \ell_{(i)\sigma_{r+1}}$ and $Y_r:=\bigoplus_{j\neq t_r} M^{\sigma_{r+1}}_j.$ Then $b_r<a_r<c_r,b_r<d_r<c_r, a_r+d_r=b_r+c_r$ and $Y_r$ contains $M\big(\sum^{t_r-1}_{i=1} \ell_{(i)\sigma_{r+1}}\big)\oplus M\big(\sum^{t_r+1}_{i=1} \ell_{(i)\sigma_{r+1}}\big)$ as a direct summand. Further, $Y_r$ has no indecomposable direct summand $Z$ with $b_r< \ell(Z)< c_r.$ It then follows from Lemma \ref{lift-exact} that there is an ${\add}(Y_r)$-split sequence
$$0\lra M\big(\ell_{(t_r+1)\sigma_{r+1}}+\sum^{t_r-1}_{i=1} \ell_{(i)\sigma_{r+1}}\big)\lra M\big(\sum^{t_r-1}_{i=1} \ell_{(i)\sigma_{r+1}}\big)\oplus M\big(\sum^{t_r+1}_{i=1} \ell_{(i)\sigma_{r+1}}\big)\lra M\big(\sum^{t_r}_{i=1} \ell_{(i)\sigma_{r+1}}\big)\lra 0.$$ Clearly, $\bigoplus^s_{j=1} M^{\sigma_r}_j= Y_r\oplus M\big(\ell_{(t_r+1)\sigma_{r+1}}+\sum^{t_r-1}_{i=1} \ell_{(i)\sigma_{r+1}}\big)$ and $\bigoplus^s_{j=1} M^{\sigma_{r+1}}_j=Y_r\oplus M\big(\sum^{t_r}_{i=1} \ell_{(i)\sigma_{r+1}}\big).$ By Lemma \ref{split-thm}, ${\End}_A(\bigoplus^s_{j=1} M^{\sigma_r}_j)$ and ${\End}_A(\bigoplus^s_{j=1} M^{\sigma_{r+1}}_j)$ are derived equivalent. $\square$

\begin{Rem}\label{mult-eq}
The sums ``$\sum^j_{i=1} \ell _i$'' and ``$\sum^j_{i=1} \ell_{(i)\sigma}$'' appearing in Lemma \ref{de} are related to the definition of $D$-equivalences of matrices. For $s\geq 2$ and a series of integers $m_s>m_{s-1}>\cdots>m_1\geq 1$, let $\ell_1:=m_1$ and $\ell_i:=m_i-m_{i-1}$ for $2\leq i\leq s$. Then $m_j=\sum^j_{i=1}\ell_i$ for $j\in [s].$ For another series of integers $n_s>n_{s-1}>\cdots>n_1\geq 1$, if  $\{\{m_s-m_{s-1},\cdots,m_1\}\}=\{\{n_s-n_{s-1},\cdots,n_1\}\}$, then there exists some $\sigma\in \Sigma_s$ such that $n_j=\sum^j_{i=1}\ell_{(i)\sigma}$ for $j\in [s].$ Moreover, if $\{\{m_s-m_{s-1},\cdots,m_1\}\}=\{\{n_s-n_{s-1},\cdots,n_1\}\}$ and if there are two irreducible polynomials $f(x)$ and $g(x)$ in $R[x]$ such that $R[x]/(f(x)^{m_s})\simeq R[x]/(g(x)^{n_s})$ as algebras, then it follows from Lemma \ref{de} that $\End_{R[x]/(f(x)^{m_s})}(\bigoplus_{k\in [s]} R[x]/(f(x)^{m_k}))$ and $\End_{R[x]/(g(x)^{n_s})}(\bigoplus_{k\in [s]} R[x]/(g(x)^{n_k}))$ are derived equivalent.
\end{Rem}

Two algebras $\Lambda$ and $\Gamma$ are said to be stably equivalent if their stable module categories $\Lambda\stmc$ and $\Gamma\stmc$ are equivalent as $R$-linear categories. Let $F$ be a stable equivalence
between  $\Lambda$ and $\Gamma$. Then $F$ induces a one-to-one correspondence between $\Lambda\modcat_{\mathscr{P}}$ and  $\Gamma\modcat_{\mathscr{P}}.$

Now, suppose that $G:\stmc{A}  \to A\stmc$ is a stable equivalence. For $n\ge 2$, $\Gamma_{n-1}:=\{M(i)\mid i\in [n-1]\}\subseteq A\modcat_{\mathscr{P}}$. Then $G$ induces an action $\ol{G}$ on $\Gamma_{n-1}$, namely, for  $M\in \Gamma_{n-1}$, $\ol{G}(M)$ is the unique module in $\Gamma_{n-1}$ such that $\ol{G}(M)\simeq G(M)$ in $A\modcat$. Clearly, $\overline{\Omega_A}(M(i))= M(n-i)$, where $\Omega_A$ is the syzygy operator of $A$.

\begin{Lem}\label{self} Let $n\ge 2$. If $G$ is a stable equivalence from $A$ to itself, then the induced action $\ol{G}$ on $\Gamma_{n-1}$ is either $\overline{\Omega_A}$ or the identity action.
\end{Lem}

{\it Proof.} If $n=2$, the conclusion is clear. Let $n\geq 3$. Since $A$ is a local, symmetric and Nakayama algebra, there are  almost split sequences
in $A\modcat$:
$$0\lra M(1)\lra M(2)\lra M(1)\lra 0$$ and
$$0\lra M(j)\lra M(j-1)\oplus M(j+1)\lra M(j)\lra 0$$ for $2\leq j\leq n-1.$ Let Irr$(X,Y)$ denote the $R$-space $\rad_A(X,Y)/\rad_A^2(X,Y)$. By a general result on stable equivalences (see \cite[Lemma 1.2, p. 336]{ARS}), we have Irr$(X,Y)\simeq$Irr$(G(X),G(Y))$ as $R$-spaces for $X,Y\in A\modcat_{\mathscr{P}}$ . It then follows that $G(M(1))\simeq M(1)$ or $G(M(1))\simeq M(n-1)=\Omega_A(M(1))$.  If $G(M(1))\simeq M(1)$, we can show that $G(M(i))\simeq M(i)$ for $i\in [n-1]$. If $G(M(1))\simeq M(n-1)=\Omega_A(M(1))$, then we consider the stable equivalence $\Omega_A^{-1}G$. Since $\Omega_A:\stmc{A}\to \stmc{A}$ is a stable equivalence, it follows from $(\Omega_A^{-1}G)(M(1))\simeq M(1)$ that $(\Omega^{-1}_AG)(M(i))\simeq M(i)$ for $i\in [n-1]$. Hence $\ol{G}$ is the identity map or equals $\overline{\Omega_A}$. $\square$

\medskip
A non-projective, non-injective simple module over an Artin algebra is called a {\it node} if the middle term of the almost split sequence starting at the simple module is projective (see \cite{MV1}).

By \cite[Lemma 1]{MV1}, a non-injective simple module $S$ of an Artin algebra is a node if and only if $S$ is not a composition factor of $\rad(Q)/\soc(Q)$ for any indecomposable projective module $Q$. Thus an Artin algebra has no nodes if and only if every non-projective, non-injective simple module is a composition factor of ${\rm rad}(P)/{\rm soc}(P)$ for some indecomposable projective module $P$.

Given an Artin algebra $\Lambda$, let $I$ be the trace of the direct sum of all non-isomorphic nodes in $\Lambda$, and $J$ be the left annihilator of $I$ in $\Lambda.$ Mart\'{i}nez-Villa showed in \cite[Theorem 2.10]{MV1} that an Artin algebra $\Lambda$ with nodes is stably equivalent to the triangular matrix algebra
$$\Lambda':=\left(
      \begin{array}{cc}    \Lambda/I & 0 \\     I & \Lambda/J \\   \end{array} \right)$$
without nodes. It is shown that $\Lambda$ and $\Lambda'$ have the same numbers of non-isomorphic, non-projective simples (see \cite[Lemma 2.10 (3)]{xz3}). We often say  that $\Lambda'$ is obtained from $\Lambda$ by eliminating nodes.

\begin{Lem}\label{node} Let $M$ be a generator for $A\modcat$.

$(1)$ If $n\geq 2$, then every simple $\End_A(M)$-module is neither projective nor injective.

$(2)$ $\End_A(M)$ has nodes if and only if $n = 2$.
\end{Lem}

{\it Proof.} Recall that $A:=R[x]/(f(x)^n)$ with $f(x)$ an irreducible polynomial in $R[x]$. Set $E:=\End_A(M)$. Since $A$ is a local, symmetric and Nakayama algebra, any indecomposable direct summand of $M$ is isomorphic to a submodule of $_AA$. Let $e_M$ denote the Hom-functor $\Hom_A(M,-): \add(_AM)\to \prj{E}$ and $P:=e_M(_AA).$ Clearly, $e_M$ is an equivalence and $P$ is an indecomposable projective-injective $E$-module with $\soc(P)\simeq \top(P)$ since $A$ is a local symmetric Nakayama algebra. The left exactness of $e_M$ implies that any indecomposable projective $E$-module is isomorphic to a submodule of $P$.

(1) Suppose $n\geq 2$. Then  the indecomposable, projective-injective $E$-module $P$ is not simple. To show that every indecomposable projective $E$-module is not simple, we show that $\Hom_E(X,Y)\ne 0$ for any indecomposable projective $E$-modules $X$ and $Y$. In fact, let $X'$ and $Y'$ be indecomposable direct summands of $_AM$ such that $e_M(X')=X$ and $e_M(Y')=Y$. Since $A$ is a local algebra, we  always have $\Hom_A(X',Y')\ne 0$. Thus $$\Hom_E(X,Y)=\Hom_{E}(e_M(X'),e_M(Y'))\simeq \Hom_A(X',Y')\ne 0.$$
This implies that $E$ has no simple projective modules. Note that the Nakayama functor $\nu:\prj{E}\to E\inj$ is an equivalence. Hence, for any indecomposable injective $E$-modules $U, V$, we also have $$\Hom_E(U,V)\simeq \Hom_E(\nu^{-1}(U),\nu^{-1}(V))\ne 0.$$ This shows that $E$ has no simple injective modules.

(2) If $n=1$, then the algebra $A$ is simple and therefore the algebra $E$ is semisimple. Thus $E$ has no nodes by definition.

If $n=2$, then $\End_A(M)$ is Morita equivalent to either $A$ or the Auslander algebra $C$ of $A$. Note that $C$ is a Nakayama algebra with $2$ indecomposable projective $C$-modules $P_1$ and $P_2=P$ of lengths $2$ and $3$, respectively. There is an almost split sequence $0\to \top(P_2)\to P_1\to \top(P_1)\to 0$, which shows that $\End_A(M)$ has a node.

Finally, we consider $n\geq 3$. Let $S$ be a non-projective, non-injective simple $E$-module, we show that $S$ is a composition factor of ${\rm rad}(P)/{\rm soc}(P).$ Actually, let $Q$ be the projective cover of $S$ with $Q=\Hom_A(M,Y)$ for $Y\in \add(M)$. Then $Q$ is a submodule of $P$ with $\soc(Q)=\soc(P)$. By (1), $Q$ is not simple. If $Q\neq P$, then $Q\subset \rad(P)$ and $0\ne Q/\soc(Q)$ is a submodule of $\rad(P)/\soc(P)$. Thus $S$ is a composition factor of $\rad(P)/\soc(P)$. If $Q=P$, then the multiplicity of $\top(P)$ in $P$ is at least $LL(A)\ge 3$, thus $S=\top(P)$ is a composition factor of $\rad(P)/\soc(P)$.  Hence $E$ has no nodes.  $\square$

\begin{Rem}\label{rmk2.11} If $n=2$ and $C$ is the Auslander algebra of $A:=R[x]/(f(x)^2)$, then $C$ has $2$ indecomposable projective modules  $P_1$ and $P_2$, the non-projective indecomposable $C$-modules are $S_1: = \top(P_1)$, $S_2:$ = top$(P_2)$ and the injective envelope $I(S_1)$ of $S_1$. Moreover, there are only two almost split sequences of $C$-modules $$0\lra S_1\lra I(S_1)\lra S_2\lra 0,\quad\quad0\lra S_2\lra P_2\lra I(S_1)\lra 0.$$ By eliminating nodes, $A$ and $C$ are stably equivalent to the corresponding triangular matrix algebras $A'$ and $C'$, respectively. The only non-injective, indecomposable projective $A'$-module is simple, while $C'$ has two non-injective, indecomposable projective modules one of which is simple. Note that the Frobenius parts of both $A'$ and $C'$ are zero.
\end{Rem}

\begin{Lem}\label{add}
Given $M,N\in A\modcat$, the endomorphism algebras ${\End}_A(M)$ and $\End_A(N)$ are Morita equivalent if and only if the basic modules $\mathcal{B}(M)$ and $\mathcal{B}(N)$ of $M$ and $N$ are isomorphic.
\end{Lem}

{\it Proof.} Suppose that ${\End}_A(M)$ and $\End_A(N)$ are Morita equivalent. Then there is an $R$-linear equivalence $G: {\add}(M)\ra {\add}(N)$. For $j\in [n]$, $\End_A(R[x]/(f(x)^j)\simeq R[x]/(f(x)^j)$ as algebras. Thus for indecomposable $A$-modules $X$ and $Y$, $\End_A(X)\simeq \End_A(Y)$ if and only if $X\simeq Y.$ It then follows from $\End_A(C)\simeq \End_A(G(C))$ for  $C\in {\add}(M)$ that $X\simeq G(X)$ for any indecomposable module $X\in \add(M)$. Therefore $\mathcal{B}(M)\simeq \mathcal{B}(N)$ as $A$-modules. $\square$

\smallskip
An algebra is said to be \emph{representation-finite} if it has only finitely many non-isomorphic indecomposable modules. Consequently, if $\add(N)\subseteq \add(M)$ and $\End_A(M)$ is representation-finite, then $\End_A(N)$ is representation-finite.

If $R$ is perfect or $f(x)\le x^r-1$  for some $r\in \mathbb{Z}_{>0}$ (for instance, $f(x)\le m_d(x)$ for an invertible matrix $d\in M_m(R)$ of finite order) , then $f(x)$ is a separable polynomial in $R[x]$.

\begin{Lem} \label{poi}
If the polynomial $f(x)$ is separable and $K:=R[x]/(f(x))$, then $A$ can be viewed as a $K$-algebra and $A\simeq K[x]/(x^n)$ as algebras over $K$.
\end{Lem}

{\it Proof.} Since $f(x)$ is separable and $\rad(A)=(f(x))/(f(x)^n)$, we know that $A/\rad(A)\simeq K$ is a separable $R$-algebra. By Wedderburn-Malcev Theorem \cite[Theorems 24 and 28]{We}, there exists a subalgebra $S$ of $A$ such that $A=S\oplus \rad(A)$ as $R$-vector spaces. Consequently, $S\simeq A/\rad(A)\simeq K.$ So $A$ can be viewed as a $K$-algebra. Since $A$ is a finite-dimensional, elementary, local $K$-algebra of representation-finite type, there is a natural number $m$ such that $A\simeq K[x]/(x^m)$. By comparing the $K$-dimensions of the algebras in this isomorphism,  we get $m=n$. $\square$

\begin{Koro} \label{St-i}
If the polynomial $f(x)$ is separable and $g(x)\in R[x]$ is irreducible such that $A$ is stably equivalent to $R[x]/(g(x)^m)$ for an integer $m\ge 2$, then $A\simeq R[x]/(g(x)^m)$.
\end{Koro}
{\it Proof.}  Since stably equivalent algebras of representation-finite type have the same number of non-projective, indecomposable modules, we get $n-1=m-1$, and therefore $n=m$. Set $B:=R[x]/(g(x)^m)$. Let $F: A\stmc{}\ra B\stmc$ be a stable equivalence and $S$ the only simple $A$-module. Thanks to $n=m\ge 2$, $S$ is not projective and $\End_A(S)\simeq \underline{\End}_A(S)$. Thus $F(S)$ is indecomposable and $$\End_A(S)\simeq \underline{\End}_A(S)\simeq \underline{\End}_B(F(S))=\End_B(F(S))/\mathcal{P}(F(S),F(S))$$is a division ring.
Since $\mathcal{P}(F(S),F(S))\subset \rad(\End_B(F(S)))$, it follows that $\mathcal{P}(F(S),F(S))=\rad(\End_B(F(S)))$. This yields the following isomorphism of algebras $$R[x]/(f(x))\simeq A/\rad(A)\simeq\underline{\End}_{A}(S)\simeq \underline{\End}_B(F(S))\simeq B/\rad(B)\simeq R[x]/(g(x)).$$
In particular, $g(x)$ is also a separable polynomial. Let $K:=R[x]/(f(x))$. Then Lemma \ref{poi} implies that $A\simeq K[x]/(x^n)\simeq B$ as algebras. $\square$

\medskip
For $c\in M_n(R)$, let $R[c]$ be the unitary subalgebra of $M_n(R)$ generated by $c$. Then there is a surjective homomorphism $\varphi: R[x]\ra R[c]$, sending $x$ to $c$. Let $R^n$ be the $n$-dimensional vector space over $R$ consisting of column vectors. Then $R^n$ is naturally an $R[c]$-module and can be viewed as an $R[x]$-module via $\varphi$. By the classification of modules over principal ideal domains, there exist irreducible polynomials $f_1(x),\cdots,f_s(x)\in R[x]$ and positive integers $t_{ij}$ such that
$$(\star)\quad \quad R^n\simeq \bigoplus^s_{i=1}\bigoplus^{l_i}_{j=1} \; R[x]/(f_i(x)^{t_{ij}})$$
as $R[x]$-modules.

\begin{Def}
The polynomials $f_i(x)^{t_{ij}}$ in $(\star)$ are called the \emph{elementary divisors} of $c$ (over $R$).
\end{Def}

Note that $\mathcal{E}_c$ is the set of pairwise distinct elementary divisors of $c$ and $\mathcal{M}_c$ is the set of maximal elements of $\mathcal{E}_c$ with respect to ``$\le$''. In particular, $m_c(x)$ is the product of elementary divisors in $\mathcal{M}_c$ and $\Ker(\varphi)$ is exactly the ideal of $R[x]$ generated by $m_c(x)$. Let $A_c:=R[x]/(\Ker(\varphi))\simeq R[c].$

\begin{Lem} \label{bijection} There is a bijection $\pi$ from $\mathcal{E}_c$ to the set of non-isomorphic indecomposable direct summands of the $A_c$-module $R^n$, sending $h(x)$ to the $A_c$-module $R[x]/(h(x))$ for $h(x)\in \mathcal{E}_c$.
\end{Lem}

Suppose that $R$ is a field of characteristic $p\geq 0$. For a positive integer $m$, there exist unique determined $s, m'\in \mathbb{N}$ such that $m=p^s m'$ and $p\nmid m'$, we define $\nu_p(m):=s$. Here, we understand $\nu_p(m):=0$ if $p=0$.

In the following, we denote by $e_{ij}$ the matrix units in $M_n(R)$ and by $c_{\sigma}:=\sum^s_{i=1} e_{i,(i)\sigma}$ the permutation matrix of $\sigma\in \Sigma_n$.

\begin{Lem}\label{per}
Let $R$ be a field of characteristic $p\ge 0$ and $\sigma\in \Sigma_n$ a permutation of cycle type $(\lambda_1,\cdots,\lambda_k)$. Then $\mathcal{E}_{c_{\sigma}} = \{f(x)^{p^{\nu_p(\lambda_i)}}\mid i\in[k], f(x)~\mbox{is irreducible}, f(x)\le x^{\lambda_i}-1\}$.
\end{Lem}

{\it Proof.} For conjugate permutations in $\Sigma_n$, their corresponding permutation matrices are similar, and therefore have the same elementary divisors. Thus, without loss of generality, we may assume that $\sigma=(1,\cdots,\lambda_1)(\lambda_1+1,\cdots,\lambda_1+\lambda_2)\cdots(\sum^{k-1}_{j=1}\lambda_j+1,\cdots,n)$. Then $c_\sigma$ is a diagonal block-matrix of the form $c_\sigma=c_{\sigma_1}\oplus c_{\sigma_2}\oplus \cdots \oplus c_{\sigma_k}$, where $\sigma_i$ is a $\lambda_i$-cycle in $\Sigma_{\lambda_i}$ for $i\in [k]$. In particular, $\mathcal{E}_{c_{\sigma}} =\bigcup_{i\in [k]}\mathcal{E}_{c_{\sigma_i}}$. For a matrix $d\in M_m(R)$, let $p_d(x)$ denote the characteristic polynomial of $d$ in $R[x]$. For $i\in [k]$, we write $\lambda_i=p^{\nu_p(\lambda_i)}\lambda_i'$ with $p\nmid \lambda'_i$. Then $x^{\lambda'_i}-1=\prod^{h_i}_{j=1}f_{ij}(x)$ is a product of distinct irreducible polynomials
$f_{ij}(x)$ in $R[x]$. It follows from $p_{c_{\sigma_i}}(x)=x^{\lambda_i}-1=x^{p^{\nu_p(\lambda_i)}\lambda_i'}-1 =(x^{\lambda'_i}-1)^{p^{\nu_p(\lambda_i)}}=
\prod^{h_i}_{j=1}f_{ij}(x)^{p^{\nu_p(\lambda_i)}}$ that
$p_{c_{\sigma_i}}(x)=m_{c_{\sigma_i}}(x)=x^{\lambda_i}-1$. Hence $\mathcal{E}_{c_{\sigma_i}}=\mathcal{M}_{c_{\sigma_i}}$. This implies $\mathcal{E}_{c_{\sigma}} = \{f(x)^{p^{\nu_p(\lambda_i)}}\mid i\in[k], f(x)~\mbox{is irreducible}, f(x)\le x^{\lambda_i}-1\}$. $\square$

\smallskip
Now, we prove a result on congruences of matrices that appear as the Cartan matrices of the endomorphism rings of modules over polynomial algebras.
Two multisets $\{\{x_1,\cdots,x_s\}\}$ and $\{\{y_1,\cdots,y_s\}\}$ are equal if and only if there exists a $\sigma\in \Sigma_s$ such that $(y_1,\cdots,y_s)^\sigma:=(y_{(1)\sigma},\cdots,y_{(s)\sigma})=(x_1,\cdots,x_s).$

\begin{Lem}\label{mat}
For an integer $s\geq 2$, let $m_1> m_2> \cdots >m_s\ge 1$ and $n_1>n_2>\cdots>n_s\ge 1$ be two series of integers with $m_1=n_1$. Set $X:=\sum^s_{k=1} (\sum^k_{l=1} m_k(e_{kl}+e_{lk})-m_k e_{kk})\in M_s(\mathbb{Z})$ and $Y:=\sum^s_{k=1} (\sum^k_{l=1} n_k(e_{kl}+e_{lk})-n_k e_{kk})\in M_s(\mathbb{Z})$. Then $X$ and $Y$ are congruent in $M_s(\mathbb{Z})$ if and only if there is $\sigma\in \Sigma_s$ such that $(n_1-n_2,\cdots,n_{s-1}-n_s,n_s)=(m_1-m_2,\cdots,m_{s-1}-m_s,m_s)^\sigma$.
\end{Lem}

{\it Proof.}  We define $U:=I-\sum^{s-1}_{t=1}e_{t,t+1}, D_1:={\rm diag}(m_1-m_2,\cdots,m_{s-1}-m_s,m_s)$ and $D_2 :={\rm diag}(n_1-n_2,\cdots,n_{s-1}-n_s,n_s).$ Then $U^{tr}XU = D_1$ and $U^{tr} Y U = D_2$ Thus $X$ and
$Y$ are congruent in $M_s(\mathbb{Z})$ if and only if $D_1$ and $D_2$ are congruent in $M_s(\mathbb{Z}).$ Now, we show that $D_1$ and $D_2$ are congruent in $M_s(\mathbb{Z})$ if and only if there is an element $\sigma\in \Sigma_s$ such that $(n_1-n_2,\cdots,n_{s-1}-n_s,n_s)=(m_1-m_2,\cdots,m_{s-1}-m_s,m_s)^\sigma.$ Indeed, if $(n_1-n_2,\cdots,n_{s-1}-n_s,n_s)=(m_1-m_2,\cdots,m_{s-1}-m_s,m_s)^\sigma,$ then $c^{tr}_{\sigma} D_1 c_{\sigma}=D_2$. This means that $D_1$ and $D_2$ are congruent in $M_s(\mathbb{Z}).$ Conversely, suppose that $D_1$ and $D_2$ are congruent in $M_s(\mathbb{Z})$. Then there is an invertible matrix $H=(a_{ij})_{1\leq i,j\leq s}\in M_n(\mathbb{Z})$ such that $H^{tr} D_1 H=D_2.$ This implies
$$(*)\quad \sum^{s-1}_{r=1} (\sum^s_{k=1}(a^2_{kr}))(m_r-m_{r+1})+(\sum^s_{k=1}(a^2_{ks}))m_s=n_1=m_1.$$
Since $H$ is invertible in $M_s(\mathbb{Z})$, each column of $H$ has a nonzero element, and therefore $\sum^s_{k=1}(a^2_{kr})\ge 1$ for $r\in [s].$ Now it follows from ($*$) that  $\sum^s_{k=1}(a^2_{kr})=1$ for all $r\in [s]$. Thus each row and column of $H$ has only one nonzero entry which is either $1$ or $-1$. This implies that $H = \epsilon c_{\tau}$ for $\tau\in \Sigma_s$ and $\epsilon$ a diagonal matrix with the entries in $\{1,-1\}$.  Hence $H^{tr}=H^{-1}.$ This shows that $D_1$ and $D_2$ are similar. So $\{\{m_1-m_2,\cdots,m_{s-1}-m_s,m_s\}\}= \{\{n_1-n_2,\cdots,n_{s-1}-n_s,n_s\}\}$ as multisets, and therefore $(n_1-n_2,\cdots,n_{s-1}-n_s,n_s)=(m_1-m_2,\cdots,m_{s-1}-m_s,m_s)^\sigma$ for some $\sigma\in \Sigma_s.$ $\square$

\subsection{Centralizer matrix algebras of representation-finite type}

We first mention a few properties of centralizer matrix algebras.

\begin{Lem}\label{iso-pr}
For $c\in M_n(R)$, the followings hold.

$(1)$ There are isomorphisms of $R$-algebras: $$S_n(c,R)\simeq S_n(c^{tr},R)\simeq S_n(c,R)^{\opp}\simeq \End_{A_c}(R^n),$$where $c^{tr}$ denotes the transpose of the matrix $c$.

$(2)$ Let $p_c(x)$ be the characteristic polynomial of $c$. Then $S_n(c,R)=R[c]$ if and only if $m_c(x)=p_c(x)$.
\end{Lem}

{\it Proof.} (1) The first isomorphism follows from the fact that any matrix over a field is similar to its transpose \cite[Theorem 66, p.76]{kaplansky}, the second isomorphism is given by sending a matrix in $S_n(c^{tr},R)$ to its transpose in $S_n(c,R)^{\opp}$, and the last isomorphism follows by interpreting $c$ as a linear transformation on the $n$-dimensional $R$-space $R^n$ with respect to a basis.

(2) Note that $m_c(x)$ is the product of polynomials in $\mathcal{M}_c$ and that $p_c(x)$ is the product of all elementary divisors (counting multiplicity) of $c$. Thus $m_c(x)=p_c(x)$ if and only if $\mathcal{M}_c$ consists of all elementary divisors of $c$. Now, suppose $m_c(x)=p_c(x)$. Then $\mathcal{M}_c=\mathcal{E}_c$ and it follows from the decomposition of $R^n$ in $(\star)$ that $R^n\simeq R[c]$ as $R[c]$-modules. Thus $S_n(c,R)\simeq \End_{R[c]}(R^n)\simeq \End_{R[c]}(R[c])\simeq R[c]$, and therefore the subalgebra $R[c]$ of $S_n(c,R)$ has the same $R$-dimension as $S_n(c,R)$ does. Hence $R[c] = S_n(c,R)$.

Suppose $S_n(c,R)=R[c]$. Then $S_n(c,R)=R[c]\simeq R[x]/(m_c(x))$. This shows that $S_n(c,R)$ is a symmetric algebra.  Hence $R^n$ is a projective generator for $R[c]\modcat$ (see, for example, the argument of the proof of Lemma \ref{lem3.1}(2)). Since $R[c]\simeq R[x]/(m_c(x))\simeq \prod_{f(x)\in \mathcal{M}_c}R[x]/(f(x))$ is a basic algebra, we have $R^n\simeq R[c]\oplus M$ for a projective $R[c]$-module $M$. Thus $R[c]=S_n(c,R)\simeq\End_{R[c]}(R^n)\simeq \End_{R[c]}(R[c]\oplus M)$. This implies $M=0$, that is, $R^n\simeq R[c]$ as $R[c]$-modules. In particular, ${\rm deg}(m_c(x))={\rm dim}_R(R[c])={\rm dim}_R(R^n)=n$. Then ${\rm deg}(p_c(x))=n= {\rm deg}(m_c(x))$. But $m_c(x)$ is a divisor of $p_c(x)$, and therefore $m_c(x)=p_c(x)$.
$\square$

\medskip
In general, $S_n(c,R)$ has neither to equal $R[c]$, nor to be representation-finite. Next, we point out a condition for $S_n(c,R)$ to be representation-finite in terms of elementary divisors.

\begin{Lem}\label{rep-f}
Suppose $R$ is a perfect field. For $c\in M_n(R)$ and $f(x)\in \mathcal{M}_c$, let $b_{f(x)}:= \mbox{max}\{3, P_c(f(x))\}$. Then $S_n(c,R)$ is representation-finite if and only if $P_c(f(x))\subseteq \{1,b_{f(x)}-1,b_{f(x)}\}$ for all $f(x)\in \mathcal{M}_c$.
\end{Lem}

{\it Proof.} Since $R$ is perfect, all irreducible factors of $m_c(x)$ are separable over $R$. Clearly, $S_n(c,R)$ is representation-finite if and only if every block of $S_n(c,R)$ is representation-finite. The blocks of $S_n(c,R)$ are parameterized by $\mathcal{M}_c$. Let $g(x)\in R[x]$ be an irreducible polynomial such that $g(x)^s\in \mathcal{M}_c$ with $s\in\mathbb{N}$. Then $b_{g(x)^s}=\max\{3,s\}$. By Lemmas \ref{M}, \ref{bijection} and \ref{iso-pr}, the block of $S_n(c,R)$ related to $g(x)^s$ is Morita equivalent to the algebra
$$E:=\End_{R[x]/(g(x)^s)}\big(\bigoplus_{t\in P_c(g(x)^s)}R[x]/(g(x)^t)\big).$$ Since $R$ is a perfect field, the algebraic closure  $\overline{R}$ of $R$ is a separable extension of $R$. By \cite[Theorem 3.3]{JL} which says that, for a separable extension $L/R$, a finite-dimensional $R$-algebra $A$ is representation-finite if and only if so is the $L$-algebra $L\otimes_R A$.  Hence it suffices to consider when $\overline{R}\otimes _R E$ is representation-finite. Since $R$ is perfect, $g(x)$ has only simple roots in $\overline{R}$. Let $\alpha_1,\cdots,\alpha_m$ be the roots of $g(x)$ in $\overline{R}$. Then
{\small $$\overline{R}\otimes _R E\simeq \End_{\overline{R}\otimes_R R[x]/(g(x)^s)}\big(\overline{R}\otimes_R \bigoplus_{t\in P_c(g(x)^s)}R[x]/(g(x)^t)\big)\simeq \End_{\overline{R}[x]/(\prod^m_{i=1}(x-\alpha_i)^s)}\big(\bigoplus_{t\in P_c(g(x)^s)}\overline{R}[x]/(\prod^m_{i=1}(x-\alpha_i)^t)\big).$$}Thus each block of $\overline{R}\otimes _R E$ is isomorphic to $\End_{\overline{R}[x]/(x^s)}\big(\bigoplus_{t\in P_c(g(x)^s)}\overline{R}[x]/(x^t)\big). $ By \cite[Theorem 2.1 (i)]{YV} (see also \cite{dr}), the endomorphism algebra $\End_{\overline{R}[x]/(x^s)}\big(\bigoplus_{t\in P_c(g(x)^s)}\overline{R}[x]/(x^t)\big) $ is representation-finite if and only if either $s\le 3$ and $P_c(g(x)^s)\subseteq \{1,2,3\}$ or $s\geq 4$ and $P_c(g(x)^s)\subseteq \{1,s-1,s\}$. This is equivalent to saying that $P_c(g(x)^s)\subseteq \{1,b_{g(x)^s}-1,b_{g(x)^s}\}$. $\square$

\medskip
As a corollary of Lemma \ref{rep-f}, we have the following.

\begin{Koro}\label{rep-p}
Let $R$ be a perfect field of characteristic $p$ and $c\in M_n(R)$ be a permutation matrix such that the associated permutation of $c$ is of the cycle type $\lambda=(\lambda_1,\cdots,\lambda_s)$. Then $S_n(c,R)$ is representation-finite if and only if there exists an positive integer $t$ such that $\nu_p(\lambda_i)\in \{0,t\}$ for all $i\in [s]$.
\end{Koro}

{\it Proof.} By Lemma \ref{per}, for $g(x)\in \mathcal{M}_c$, all the integers in $P_c(g(x))$ are $p$-powers and $(x-1)^{p^{\nu_p(\lambda_i)}}$ is an elementary divisor of $c$ for $i\in [s]$. In particular, $(x-1)^{p^m}\in \mathcal{M}_c$ with $m=\max\{\nu_p(\lambda_i)\mid i\in [s]\}$.

Now, suppose that $S_n(c,R)$ is representation-finite. By Lemma \ref{rep-f}, we deduce that $P_c((x-1)^{p^m})$ does not contain two different $p$-powers $p^a>1$ and $p^b>1$. Since $p^{\nu_p(\lambda_i)}\in P_c((x-1)^{p^m})$ for $i\in [s]$, there do not exist $\lambda_i$ and $\lambda_j$ with $i,j\in [s]$ such that $\nu_p(\lambda_i)>\nu_p(\lambda_j)\geq 1$, that is, there exists an integer $t>0$ such that $\nu_p(\lambda_i)\in \{0,t\}$ for all $i\in [s]$.

Conversely, suppose that there exists an integer $t>0$ such that $\nu_p(\lambda_i)\in \{0,t\}$ for all $i\in [s]$. Then, for each $g(x)\in \mathcal{M}_c$, we deduce from Lemma \ref{per} that $P_c(g(x))\subseteq \{1,p^t\}$. Thus it follows from Lemma \ref{rep-f} that $S_n(c,R)$ is representation-finite. $\square$

\section{Derived and stable equivalences of centralizer matrix algebras\label{Pf}}
This section is devoted to proving all results mentioned in the introduction.

Assume that $R$ is a field of characteristic $p\geq 0$ unless stated otherwise. For $c\in M_n(R)$,  let $m_{R,c}(x)$ or $m_c(x)$ be the minimal polynomial of $c$ over $R$ and $A_c:=R[x]/(m_c(x))$. Now, let $d\in M_m(R)$. We write
$$m_c(x)=\prod^{l_c}_{i=1} f_i(x)^{n_i} \mbox{  for } n_i\ge 1 \; \mbox{ and } \; m_d(x)=\prod^{l_d}_{j=1} g_j(x)^{m_j} \; \mbox{  for } m_j\ge 1,$$   $$U_i:=R[x]/(f_i(x)^{n_i})\mbox{  for } i\in [l_c] \; \mbox{ and }\; \; V_j:=R[x]/(g_j(x)^{m_j})\mbox{ for } j\in [l_d],$$ where $f_i(x)$ and $g_j(x)$ are irreducible polynomials in $R[x]$. Thus $U_i$ and $V_j$ are local, symmetric Nakayama $R$-algebras, and
$$A_c\simeq U_1\times U_2\times\cdots\times U_{l_c} \mbox{ and }\; A_d \simeq V_1\times V_2\times\cdots\times V_{l_d}.$$

Recall that $A_c\simeq R[c]$ and $R^n=\{(a_1,a_2,\cdots,a_n)^{tr}\mid a_i\in R, 1\le i\le n\}$ is viewed as an $A_c$-module. According to these blocks of $A_c$ and $A_d$, we decompose the $A_c$-module $R^n$ and the $A_d$-module $R^m$ as
$$R^n = \bigoplus^{l_c}_{i=1} M_i \; \mbox{ and } \; R^m =\bigoplus^{l_d}_{j=1} N_j,$$where $M_i$ (respectively, $N_j$) is the sum of indecomposable direct summands of $R^n$ (respectively, $R^m$) belonging to the block $U_i$ (respectively, $V_j$). Then $\mathcal{B}(M_i)\simeq \bigoplus_{r\in {P_c(f_i(x)^{n_i})}} R[x]/(f_i(x)^r)$ as $U_i$-modules and $\mathcal{B}(N_j)\simeq \bigoplus_{s\in {P_d(g_j(x)^{m_j})}} R[x]/(g_j(x)^s)$ as $V_j$-modules. Since $R^n$ is a faithful $M_n(R)$-module, $R^n$ is also a faithful $R[c]$-module, and therefore $M_i$ is a faithful $U_i$-module for $i\in[l_c]$. Similarly, $N_j$ is a faithful $V_j$-module for $j\in [l_d]$.
Further, we set $$A_i:=\End_{U_i}(M_i) \; \mbox{  and } \; B_j:=\End_{V_j}(N_j)$$ for $i\in [l_c]$ and $j\in [l_d]$. Then $A_i$ and $B_j$ are indecomposable as algebras for $i\in [l_c]$ and $j\in [l_d]$. Clearly,  $A_i$ (respectively, $B_j$) is semisimple if and only if $n_i=1$ (respectively, $m_j=1$). In this case, $A_i\simeq M_k(R[x]/(f_i(x)))$ for some $k\in \mathbb{N}$ (respectively, $B_j\simeq M_t(R[x]/(g_j(x)))$ for some $t\in \mathbb{N}$). By Lemma \ref{iso-pr},
$$S_n(c,R)\simeq \prod^{l_c}_{i=1}{\End}_{U_i}(M_i)= \prod^{l_c}_{i=1}A_i \, \mbox{ and } \; S_m(d,R)\simeq \prod^{l_d}_{j=1}{\End}_{V_j}(N_j)= \prod^{l_d}_{i=1}B_j.$$

As the $R[c]$-module $R^n$ is a generator, we see that the bimodule $_{R[c]}R^n_{S_n(c,R)}$ has the double centralizer property. In particular, $\End_{S_n(c,R)}(R^n_{_{S_n(c,R)}})=R[c]$.

\subsection{Characterizations of Morita and derived equivalences: Proof of Theorem \ref{main1}}
In this section we prove the main result, Theorem \ref{main1}.

\begin{Lem} \label{3.1} $(1)$ $\mathcal{M}_c=\{f_i(x)^{n_i}\mid 1\le i\le l_c\}$ and $|\mathcal{M}_c|=l_c$.

$(2)$ If $A_i$ and $B_j$ are derived equivalent, then $U_i\simeq V_j$ and $n_i=m_j$.
\end{Lem}

{\it Proof.} (1) follows by definition. (2) is a consequence of Lemma \ref{iso}. $\square$

\begin{Lem}\label{E}
Let $c\in M_n(R)$ and $d\in M_m(R)$. Then $c\stackrel{M}\sim d$ if and only if there is an isomorphism $\varphi:R[c]\simeq R[d]$ of algebras such that $\mathcal{B}(R^n)\simeq \mathcal{B}(R^m)$, where $R^m$ is viewed as an $R[c]$-module via $\varphi$.
\end{Lem}

{\it Proof.} Suppose $c\stackrel{M}\sim d$. Then, by definition, there is a bijection $\pi$ between $\mathcal{M}_c$ and $\mathcal{M}_d$ such that $R[x]/(f(x))\simeq R[x]/((f(x))\pi)$ as algebras and ${P_c(f(x))}= {P_d((f(x))\pi)}$ for all $f(x)\in \mathcal{M}_c.$ It follows from $$R[c]\simeq \prod_{f(x)\in \mathcal{M}_c} R[x]/(f(x)) \mbox{  and }\; R[d]\simeq \prod_{g(x)\in \mathcal{M}_d} R[x]/(g(x))$$that there is an isomorphism $\varphi: R[c]\simeq R[d].$ After reordering the factors, we may assume that $(f_i(x)^{n_i})\pi=g_i(x)^{m_i}$ for $i\in [l_c].$ Then the condition ${P_c(f(x))}= {P_d((f(x))\pi)}$ implies that $\mathcal{B}(M_i)\simeq \mathcal{B}(N_i)$  for $i\in [l_c]$. Here $N_i$ is viewed as an $R[c]$-module via $\varphi$. Hence $\mathcal{B}(R^n)\simeq \mathcal{B}(R^m)$, where $R^m$ is viewed as an $R[c]$-module via $\varphi$.

Conversely, suppose that there is an isomorphism $\varphi:R[c]\simeq R[d]$ such that $\mathcal{B}(R^n)\simeq \mathcal{B}(R^m)$ when $R^m$ is regarded as an $R[c]$-module via $\varphi$. We may assume that $\varphi$ restricts to an isomorphism $\varphi_i: U_i\simeq V_i$ for $i\in [l_c].$ Then the condition
$\mathcal{B}(R^n)\simeq \mathcal{B}(R^m)$ implies that $\mathcal{B}(M_i)\simeq \mathcal{B}(N_i)$ for $i\in [l_c].$ Since $\mathcal{B}(M_i)\simeq \bigoplus_{r\in {P_c(f_i(x)^{n_i})}} R[x]/(f_i(x)^r)$ as $U_i$-modules and $\mathcal{B}(N_i)\simeq \bigoplus_{s\in {P_d(g_i(x)^{m_i})}} R[x]/(g_i(x)^s)$ as $V_i$-modules, we have ${P_c(f_i(x)^{n_i})}= {P_d(g_i(x)^{m_i})}.$  Now we define a map $\pi: \mathcal{M}_c\to \mathcal{M}_d$ by $f_i(x)^{n_i}\mapsto g_i(x)^{m_i}$ for $i\in [l_c].$ Then $\pi$ defines an $M$-equivalence $c\stackrel{M}\sim d$. $\square$

\medskip
Recall that for an algebra $A$ and an idempotent $e\in A$, the Schur functor $Ae\otimes_{eAe}-: eAe\modcat \to A\modcat$ is
fully faithful. Following \cite[Section 2]{hx3}, we say that a stable equivalence $\Phi: A\stmc{}\ra B\stmc$ of Morita type lifts to a Morita equivalence if there is a Morita equivalence $F: A\modcat{}\ra B\modcat$ such the  following diagram of functors is commutative (up to natural isomorphism)
$$\xy
(0,15)*+{\stmc{A}}="a",
(25,15)*+{\stmc{B}}="b",
(0,0)*+{A\modcat}="d",
(25,0)*+{B\modcat}="e",
{\ar^{can.} "d";"a"},
{\ar^{can.} "e";"b"},
{\ar^{\Phi}, "a";"b"},
{\ar^{F}, "d";"e"},
\endxy$$

The following is proved in \cite[Lemma 3.1, Propositions 3.3 and 3.5, and Remark 4.6]{hx3}.
\begin{Lem}\label{Sta-M}
Let $A$ and $B$ be two algebras without nonzero semisimple direct summands such that $A/\rad(A)$ and $B/\rad(B)$ are separable, and let $e$ and $f$ be $\nu$-stable idempotent in $A$ and $B$ such that $eAe$ and $fAf$ are their Frobenius parts. Suppose there is a stable equivalence $\Phi: A\stmc{}\ra B\stmc$ of Morita type. Then the followings hold.

$(1)$ If $\Phi(S)$ is isomorphic in $B\stmc$ to a simple $B$-module for each simple $A$-module $S$, then $\Phi$ lifts to a Morita equivalence.

$(2)$ The functor $\Phi$ restricts to a stable equivalence $\Phi_1: eAe\stmc{}\ra fBf\stmc$ of Morita type such that the following diagram is commutative (up to natural isomorphism)
$$\xy
(0,15)*+{\stmc{A}}="a",
(25,15)*+{\stmc{B}}="b",
(0,0)*+{\stmc{eAe}}="d",
(25,0)*+{\stmc{fBf}}="e",
{\ar^{\lambda} "d";"a"},
{\ar^{\lambda} "e";"b"},
{\ar^{\Phi}, "a";"b"},
{\ar^{\Phi_1}, "d";"e"},
\endxy$$
where $\lambda$ stands for the corresponding Schur functor. Moreover, if $\Phi_1$ lifts to a Morita equivalence, then so is $\Phi$.
\end{Lem}

\medskip
{\bf Proof of Theorem \ref{main1}}. If $S_n(c,R)$ and $S_m(d,R)$ are Morita (or derived, or almost $\nu$-stable derived) equivalent, then they have the same number of blocks, that is, $l_c=l_d$. Further, we may assume that $A_i$ and $B_i$ are Morita (or derived, or almost $\nu$-stable derived)) equivalent and that $F_i$ is such an equivalence for $i\in [l_c]$. By Lemma \ref{3.1}(2), there is an isomorphism $\varphi_i: U_i\simeq V_i$ of algebras and $n_i=m_i$ for $i\in [l_c].$

(1) Suppose $c\stackrel{M}\sim d$. Then it follows from Lemmas \ref{M} and \ref{E} that $S_n(c,R)$ and $S_m(d,R)$ are Morita equivalent. Conversely, suppose that $S_n(c,R)$ and $S_m(d,R)$ are Morita equivalent. Then it follows from Lemma \ref{add} that $\mathcal{B}(M_i)\simeq \mathcal{B}(N_i)$ if $N_i$ is regarded as a $U_i$-module via $\varphi_i$. Let $\psi_c: A_c\ra R[c]$ and $\psi_d: A_d\ra R[d]$ be the canonical isomorphisms of algebras. Then $\varphi:=\psi^{-1}_c(\prod^{l_c}_{i=1} \varphi_i)\psi_d$ is an isomorphism between the algebras $R[c]$ and $R[d]$. Hence $\mathcal{B}(R^n)\simeq \mathcal{B}(R^m)$ where $R^m$ is viewed as an $R[c]$-module via $\varphi$. By Lemma \ref{E}, we have $c\stackrel{M}\sim d.$

(2) Suppose $c\stackrel{D}\sim d$. By the definition of $D$-equivalences, $A_c\simeq A_d$ as algebras and there is a map $\pi:\mathcal{M}_c\to \mathcal{M}_d$ such that $\mathcal{H}_{P_c(f_i(x)^{n_i})}= \mathcal{H}_{P_d((f_i(x)^{n_i})\pi)}$ for $f_i(x)^{n_i}\in \mathcal{M}_c.$ Without loss of generality, we assume $(f_i(x)^{n_i})\pi=g_i(x)^{m_i}$ for $i\in [l_c].$ Then, for $i\in [l_c]$, $R[x]/(f_i(x)^{n_i})\simeq R[x]/(g_i(x)^{m_i})$ as algebras and $\mathcal{H}_{P_c(f_i(x)^{n_i})}= \mathcal{H}_{P_d(g_i(x)^{m_i})}$. Recall that $\mathcal{B}(M_i)\simeq \bigoplus_{r\in {P_c(f_i(x)^{n_i})}} R[x]/(f_i(x)^r)$ as $U_i$-modules and $\mathcal{B}(N_i)\simeq \bigoplus_{s\in {P_d(g_i(x)^{m_i})}} R[x]/(g_i(x)^s)$ as $V_i$-modules. It follows from Remark \ref{mult-eq} that $\End_{U_i}(\mathcal{B}(M_i))$ and $\End_{V_i}(\mathcal{B}(N_i))$ are derived equivalent. Thanks to Lemma \ref{M},  $A_i$ and $B_i$ are also derived equivalent. Thus $S_n(c,R)$ and $S_m(d,R)$ are derived equivalent.

Conversely, suppose that $S_n(c,R)$ and $S_m(d,R)$ are derived equivalent. Let $i\in [l_c]$. Then $A_i$ and $B_i$ are derived equivalent, and there is an isomorphism $U_i\simeq V_i$ of algebras such that $U_i/\rad(U_i)\simeq V_i/\rad(V_i)$, that is, $R[x]/(f_i(x))\simeq R[x]/(g_i(x)).$ Let $K_i$ be a splitting field for $f_i(x)g_i(x).$ Since $K_i\otimes_R A_i\simeq \End_{K_i\otimes_R U_i}(K_i\otimes_R M_i)$ and $K_i\otimes_R B_i\simeq \End_{K_i\otimes_R V_i}(K_i\otimes_R N_i)$, the two algebras $\End_{K_i\otimes_R U_i}(K_i\otimes_R M_i)$ and $\End_{K_i\otimes_R V_i}(K_i\otimes_R N_i)$ are derived equivalent.

For the irreducible polynomial $f_i(x)$, there is a separable irreducible polynomial $u_i(x)\in R[x]$ and an integer $s_i\in \mathbb{N}$ such that $f_i(x)=u_i(x^{p^{s_i}})$. Here, for $p=0$, we understand $p^{s_i}=1$. Similarly, there is a separable irreducible polynomial $v_i(x)$ and an integer $t_i\in \mathbb{N}$ such that $g_i(x)=v_i(x^{p^{t_i}})$. It follows from $K_i\otimes_R \big(R[x]/(f_i(x))\big)\simeq K_i\otimes_R \big(R[x]/(g_i(x))\big)$ that $s_i=t_i$ and that $u_i(x)$ and $v_i(x)$ have the same number of roots. Therefore $f_i(x), g_i(x),u_i(x)$ and $v_i(x)$ have the same number of distinct roots in $K_i$. Let $w_i$ be the number of roots of $u_i(x)$ in $K_i$. Suppose that $\alpha_{i1}, \alpha_{i2},\cdots, \alpha_{i w_i}$ are the roots of $f_i(x)$ in $K_i$ and that $\beta_{i1}, \beta_{i2}, \cdots, \beta_{i w_i}$ are the roots of $g_i(x)$ in $K_i$. Then $K_i\otimes_R U_i= K_i\otimes_R \big(R[x]/(f_i(x)^{n_i})\big)\simeq \prod^{w_i}_{q=1} K_i[x]/((x-\alpha_{iq})^{n_i\cdot p^{s_i}})$. Similarly, $K_i\otimes_R V_i=K_i\otimes_R \big(R[x]/(g_i(x)^{m_i})\big)\simeq \prod^{w_i}_{q=1} K_i[x]/((x-\beta_{iq})^{m_i\cdot p^{s_i}}).$ Now, we shall show $\mathcal{H}_{P_c(f_i(x)^{n_i})}= \mathcal{H}_{P_d(g_i(x)^{m_i})}$.
Indeed, given a $U_i$-module $R[x]/(f_i(x)^r)$, we have the isomorphism
$$K_i\otimes_R \big(R[x]/(f_i(x)^r)\big)\simeq \bigoplus^{w_i}_{q=1} K_i[x]/\big((x-\alpha_{iq})^{rp^{s_i}}\big)$$
as $\prod^{w_i}_{q=1} K_i[x]/\big((x-\alpha_{iq})^{n_ip^{s_i}}\big)$-modules.
Since $\Hom_{U_i}(M_i,-): \add(M_i)\to A_i\mbox{-proj}$ is an equivalence, we see that $|P_c(f_i(x)^{n_i})|$ equals the number of indecomposable projective $A_i$-modules, hence equals  the number of simple $A_i$-modules. Since derived equivalent algebras have the same number of simple modules, we get $|P_c(f_i(x)^{n_i})|=|P_d(g_i(x)^{n_i})|.$ Put $h_i=|P_c(f_i(x)^{n_i})|.$ For  $h_i=1$, we have $\mathcal{H}_{P_c(f_i(x)^{n_i})}= \mathcal{H}_{P_d(g_i(x)^{m_i})}$. So we assume that $h_i\geq 2$ and $P_c(f_i(x)^{n_i})=\{u_{i1},\cdots,u_{ih_i}\}$ with $u_{i1}> \cdots > u_{ih_i}.$ Since $A_i=\End_{U_i}(M_i)$ is Morita equivalent to $\End_{U_i}(\mathcal{B}(M_i))$, the algebra $K_i\otimes_R A_i$ is Morita equivalent to  $K_i\otimes_R \End_{U_i}(\mathcal{B}(M_i))\simeq  \End_{K_i\otimes_R U_i}(K_i\otimes_R \mathcal{B}(M_i))$. As $\mathcal{B}(M_i)\simeq \bigoplus_{k\in [h_i]} R[x]/(f_i(x)^{u_{ik}})$ as $U_i$-modules, we get
$$K_i\otimes_R \mathcal{B}(M_i)\simeq \bigoplus^{w_i}_{q=1}\bigoplus_{k\in [h_i]} K_i[x]/((x-\alpha_{iq})^{u_{ik}p^{s_i}})$$as $\prod^{w_i}_{q=1} K_i[x]/\big((x-\alpha_{iq})^{n_ip^{s_i}}\big)$-modules. Then a block of $K_i\otimes_R A_i$ is Morita equivalent to the algebra
$$E_{c,i}:=\End_{K_i[x]/((x-\alpha_{iq})^{n_ip^{s_i}})}\big(\bigoplus_{k\in [h_i]} K_i[x]/((x-\alpha_{iq})^{u_{ik}p^{s_i}})\big)$$
for some $q\in [w_i].$ Similarly,
$$K_i\otimes_R \big(R[x]/(g_i(x)^r)\big)\simeq \bigoplus^{w_i}_{q=1} K_i[x]/\big((x-\beta_{iq})^{rp^{s_i}}\big)$$
as $\prod^{w_i}_{q=1} K_i[x]/\big((x-\beta_{iq})^{n_ip^{s_i}}\big)$-modules. We write $P_d(g_i(x)^{n_i})=\{v_{i1},\cdots,v_{ih_i}\}$ with $v_{i1}> \cdots > v_{ih_i}$. Then a block of  $K_i\otimes_R B_i$ is Morita equivalent to an algebra of the form
$$E_{d,i}:=\End_{K_i[x]/((x-\beta_{iq'})^{n_ip^{s_i}})}\big(\bigoplus_{k\in [h_i]} K_i[x]/((x-\beta_{iq'})^{v_{ik}p^{s_i}})\big)$$
for some $q'\in [w_i].$ Remind that $u_{i1}=n_i=m_i=v_{i1}$, $E_{c,i}\simeq \End_{K_i[x]/(x^{n_ip^{s_i}})}(\oplus_{k\in [h_i]} K_i[x]/(x^{u_{ik}p^{s_i}}))$ and $E_{d,i}\simeq\End_{K_i[x]/(x^{n_ip^{s_i}})}(\oplus_{k\in [h_i]} K_i[x]/(x^{v_{ik}p^{s_i}})).$
Thus the Cartan matrices of  $E_{c,i}$  and $E_{d,i}$ (as $K_i$-algebras) are the $h_i\times h_i$ matrices
$$H_i:=p^{s_i}\sum^{h_i}_{k=1} (\sum^k_{l=1} u_{ik}(e_{kl}+e_{lk})-u_{ik} e_{kk}) \mbox{ and }  J_i:=p^{s_i}\sum^{h_i}_{k=1} (\sum^k_{l=1} v_{ik}(e_{kl}+e_{lk})-v_{ik} e_{kk}),$$ respectively. Then there exists an invertible matrix $\Phi_i\in M_{h_i}(\mathbb{Z})$ such that $\Phi^{tr}_i H_i\Phi_i=J_i.$ This follows from \cite[Chapter 6, Proposition 6.8.9]{Z} which says that  the Cartan matrices of derived equivalent, split algebras are congruent by an invertible matrix with integral entries. Thanks to Lemma \ref{mat}, we have $\mathcal{H}_{P_c(f_i(x)^{n_i})}= \mathcal{H}_{P_d(g_i(x)^{m_i})}$ as multisets.

Now let $\pi: \mathcal{M}_c\to \mathcal{M}_d$ be the map given by $f_i(x)^{n_i}\mapsto g_i(x)^{m_i}$ for $i\in [l_c].$ Then $\pi$ gives rise to a $D$-equivalence $c\stackrel{D}\sim d$.

(3) Suppose $c\stackrel{AD}\sim d.$ Then $U_i\simeq V_i$ for $i\in [l_c]$ by definition. The condition ${P_c(f_i(x)^{n_i})}= {P_d(g_i(x)^{m_i}))}$ or ${P_c(f_i(x)^{n_i})} = \mathcal{J}_{P_d(g_i(x)^{m_i}))}$ implies that either $\mathcal{B}(M_i)_{\mathscr{P}}\simeq \mathcal{B}(N_i)_{\mathscr{P}}$ or $\mathcal{B}(M_i)_{\mathscr{P}}\simeq\Omega_{V_i}(\mathcal{B}(N_i)_{\mathscr{P}})$ as $U_i$-modules. It then follows from Lemma \ref{alm} that $A_i$ and $B_i$ are almost $\nu$-stable derived equivalent. Hence $S_n(c,R)$ and $S_m(d,R)$ are almost $\nu$-stable derived equivalent.

Conversely, suppose that $S_n(c,R)$ and $S_m(d,R)$ are almost $\nu$-stable derived equivalent. Thanks to Lemma \ref{almst}, the almost $\nu$-stable derived equivalence $F_i$ induces a stable equivalence, say $\overline{F_i}$, between $U_i$ and $V_i$ such that $\overline{F_i}(\mathcal{B}(M_i)_{\mathscr{P}})\simeq \mathcal{B}(N_i)_{\mathscr{P}}$ and $m_i=n_i$ for $i\in [l_c]$. The isomorphism $\varphi_i$ induces a
natural stable equivalence $\Phi_i$ between $V_i$ and $U_i$. The composition $\Phi_i\circ\overline{F_i}$ is then a stable equivalence from $U_i$ to itself. Now, by Lemma \ref{self}, we deduce either $\mathcal{B}(M_i)_{\mathscr{P}}\simeq \mathcal{B}(N_i)_{\mathscr{P}}$ or $\mathcal{B}(M_i)_{\mathscr{P}}\simeq \Omega_{V_i}(\mathcal{B}(N_i)_{\mathscr{P}})$ as $U_i$-modules, where $N_i$ is viewed as a $U_i$-module via $\varphi_i$. Note that $\mathcal{B}(M_i)\simeq \bigoplus_{r\in {P_c(f_i(x)^{n_i})}} R[x]/(f_i(x)^r)$ as $U_i$-modules and $\mathcal{B}(N_j)\simeq \bigoplus_{s\in {P_d(g_j(x)^{m_j})}} R[x]/(g_j(x)^s)$ as $V_j$-modules. Thus $\mathcal{B}(M_i)_{\mathscr{P}}\simeq \mathcal{B}(N_i)_{\mathscr{P}}$ is equivalent to ${P_c(f_i(x)^{n_i})}= {P_d(g_i(x)^{m_i})}$, and $\mathcal{B}(M_i)_{\mathscr{P}}\simeq \Omega_{V_i}(\mathcal{B}(N_i)_{\mathscr{P}}$ is equivalent to ${P_c(f_i(x)^{n_i})}=\mathcal{J}_{P_d(g_i(x)^{m_i})}$ for $i\in [l_c]$.
Now we define a map $\pi: \mathcal{M}_c\to \mathcal{M}_d$ by $f_i(x)^{n_i}\mapsto g_i(x)^{m_i}$ for $i\in [l_c].$ Then $\pi$ defines an $AD$-equivalence: $c\stackrel{AD}\sim d$.

(4) Assume that either $R$ is perfect or both $c$ and $d$ are invertible matrices of finite order. Then all irreducible factors of $m_c(x)$ and $m_d(x)$ are separable polynomials over $R$. Let $A_i=\End_{U_i}(M_i)$ be a block in $S_n(c,R)$ and $e_{M_i}:=\Hom_{U_i}(M_i,-): U_i\modcat{}\ra A_i\modcat$ be an evaluation functor. For any indecomposable projective $A_i$-module $e_{M_i}(X)$ with $X$ an indecomposable direct summand of $M_i$, we have $$\End_{A_i}\big(\top(e_{M_i}(X))\big)\simeq\End_{A_i}(e_{M_i}(X))/\rad (\End_{A_i}(e_{M_i}(X)))\simeq \End_{U_i}(X)/\rad (\End_{U_i}(X))\simeq R[x]/(f_i(x)).$$ Thus $\End_{A_i}\big(\top(e_{M_i}(X))\big)$ is separable. This implies that the semisimple quotient $A_i/\rad(A_i)$ of $A_i$ is separable. Hence all the semisimple quotients of blocks of $S_n(c,R)$ and $S_m(d,R)$ are separable. In particular, all the semisimple blocks of $S_n(c,R)$ and $S_m(d,R)$ are separable.

Now, suppose that $c$ and $d$ are $SM$-equivalent. By Lemma \ref{alm} and the proof of (3), there is an almost $\nu$-stable derived equivalence between the non-semisimple blocks of $S_n(c,R)$ and $S_m(d,R)$ . Further, by \cite[Theorem 1.1]{hx1}, there is a stable equivalence $F$ of Morita type between the non-semisimple blocks $S_n(c,R)$ and $S_m(d,R)$. Since the semisimple blocks of $S_n(c,R)$ and $S_m(d,R)$ are separable algebras, $F$ can be extended to a stable equivalence of Morita type between $S_n(c,R)$ and $S_m(d,R)$.

Conversely, suppose that $S_n(c,R)$ and $S_m(d,R)$ are stably equivalent of Morita type, and that the stable equivalence is given by a functor $F$. Let $A_1,\cdots,A_s$ be the non-semisimple blocks of $S_n(c,R)$, and let $B_1,\cdots,B_s$ be the non-semisimple blocks of $S_m(d,R)$. By \cite[Theorem 2.2]{Liu1} and \cite[Lemma 4.8]{Liu2},  we may assume that $F$ induces a stable equivalence $F_i$  of Morita type, between $A_i$ and $B_i$ for $1\le i\le s$.

To show $c\stackrel{SM}\sim d$, we consider the generator $M_i$ for $U_i\modcat$. It follows from $\nu_{A_i}\Hom_{U_i}(M_i, U_i) \simeq
\Hom_{U_i}(M_i,\nu_{U_i}U_{i})$ (see \cite[Remark 2.9 (2)]{hx3}) that the Frobenius parts of $A_i$ and $B_i$ are Morita equivalent to $U_i$ and $V_i$,
respectively. Since $A_i/\rad(A_i)$ and $B_i/\rad(B_i)$ are separable, it follows from Lemma \ref{Sta-M}(2) that $F_i$ restricts to a stabe equivalence $G_i$ of Morita type between $U_i$ and $V_i$. As $f_i(x)$ is separable and both $A_i$ and $B_i$ are non-semisimple, Corollary \ref{St-i} implies that $U_i\simeq V_i$, that is, $R[x]/(f_i(x)^{n_i})\simeq R[x]/(g_i(x)^{m_i})$, and $n_i=m_i$.

Now we regard $V_i$-modules as $U_i$-modules via this isomorphism. Let $\overline{A_i}:=\End_{U_i}(U_i\oplus \mathcal{B}(M_i)_{\mathscr{P}})$,  $\overline{B_i}:=\End_{V_i}(V_i\oplus \mathcal{B}(N_i)_{\mathscr{P}})$ and $\overline{C_i}:=\End_{V_i}(V_i\oplus \Omega_{V_i}(\mathcal{B}(N_i)_{\mathscr{P}}))$, and let $e,f$ and $g$ be the $\nu$-stable idempotents of $\overline{A_i},\overline{B_i}$ and $\overline{C_i}$, defining their Frobenius parts, respectively. Then $A_i,\overline{A_i},B_i,\overline{B_i}$ and $\overline{C_i}$ are stably equivalent of Morita type, and there is the following commutative (up to natural isomorphism) diagram by Lemma \ref{Sta-M}(2):
$$\xy
(0,15)*+{\stmc{\overline{A_i}}}="a",
(25,15)*+{\stmc{\overline{B_i}}}="b",
(50,15)*+{\stmc{\overline{C_i}}}="c",
(0,0)*+{\stmc{e\overline{A_i}e}}="d",
(25,0)*+{\stmc{f\overline{B_i}f}}="e",
(50,0)*+{\stmc{g\overline{C_i}g}}="f",
{\ar^{\lambda} "d";"a"},
{\ar^{\lambda} "e";"b"},
{\ar^{\lambda} "f";"c"},
{\ar^{\Phi}, "a";"b"},
{\ar^{\Phi_1}, "d";"e"},
{\ar^{\Psi}, "b";"c"},
{\ar^{\Psi_1}, "e";"f"},
\endxy$$
where $\lambda$ is the full embedding of stable module categories induced by the corresponding Schur functor and where $\Phi$ and $\Psi$ define stable equivalences of Morita type between $\overline{A_i}$ and $\overline{B_i}$, and between $\overline{B_i}$ and $\overline{C_i}$, respectively, while $\Phi_1$ and $\Psi_1$ are the restrictions of $\Phi$ and $\Psi$, respectively. Note that $e\overline{A_i}e\simeq U_i\simeq V_i\simeq f\overline{B_i}f\simeq g\overline{C_i}g$. Identifying $f\overline{B_i}f$ with $g\overline{C_i}g$, we can choose $\Psi$ so that $\Psi_1$ is the syzygy functor on $f\overline{B_i}f\stmc$ (see the arguments in \cite[Proposition 3.3 and Corollary 3.4]{LX3}). Let $S$  be the simple $e\overline{A_i}e$-module. Then it follows from Lemma \ref{self} that either $\Phi_1(S)$ or $\Psi_1\circ \Phi_1(S)$ is simple. By Lemma \ref{Sta-M}(1), either $\Phi_1$ or $\Psi_1\circ \Phi_1$ can be lifted to a Morita equivalence, and therefore
either $\Phi$ or $\Psi\circ \Phi$ can be lifted to a Morita equivalence. It then follows from Lemma \ref{add} that either $\mathcal{B}(M_i)_{\mathscr{P}}\simeq \mathcal{B}(N_i)_{\mathscr{P}}$ or $\mathcal{B}(M_i)_{\mathscr{P}}\simeq \Omega_{V_i}(\mathcal{B}(N_i)_{\mathscr{P}})$.
Therefore ${P_c(f_i(x)^{n_i})}= {P_d(g_i(x)^{m_i}))}$ or ${P_c(f_i(x)^{n_i})}=\mathcal{J}_{P_d(g_i(x)^{m_i}))}.$ Now we define a map $\pi: \mathcal{R}_c\to \mathcal{R}_d$ by $f_i(x)^{n_i}\mapsto g_i(x)^{m_i}$ for $f_i(x)^{n_i}\in \mathcal{R}_c$. Then $\pi$ defines an $SM$-equivalence $c\stackrel{SM}\sim d$ of matrices. $\square$

Instead of $R$ being a field, we can prove the following for  noetherian domains.

\begin{Rem}\label{rmk3.3}
Suppose that $R$ is a noetherian domain, $c\in M_n(R)$ and $d\in M_m(R)$. If $S_n(c,R)$ and $S_m(d,R)$ are derived equivalent, then $c\stackrel{D}\sim d$ as matrices over the fraction field of $R$.
\end{Rem}

{\it Proof.} Assume that $R$ is a noetherian domain with $K$ its fractional field. Then it follows from $S_n(c,R)\subseteq M_n(R)$ that $S_n(c,R)$ is a finitely generated $R$-algebra. Thus $S_n(c,R)$ is a noetherian algebra, $S_n(c,R)\modcat$ is an abelian category and $\Db{S_n(c,R)}$ is well defined.

Regarding $K$ as an $R$-algebra, we have the isomorphism of $K$-algebras
\begin{equation}
\begin{aligned}
\varphi: K\otimes_R M_n(R)&\lra M_n(K), \; \;
\sum^s_{i=1} a_i\otimes b_i&\mapsto \, \sum^s_{i=1} (a_iI_n) b_i
\end{aligned}
\nonumber
\end{equation}
where $I_n$ is the identity matrix in $M_n(K).$ Further, $K$ is a flat $R$-module and there is the commutative diagram of $K$-algebras
$$\xymatrix{
K\otimes_RS_n(c,R) \ar[r]^-{\mu}\ar@{^{(}->}[d] & S_n(c,K)\ar@{^{(}->}[d]\\
K\otimes_R M_n(R) \ar[r]_-{\sim}^-{\varphi} & M_n(K)\\
}$$where $\mu$ is the restriction of $\varphi$. Remark that $\Img(\mu)$ belongs to $S_n(c,K)$. Since $K$ is the fractional field of $R$,  we can find an element $0\ne r\in R$ for each matrix $a\in M_n(K)$ such that $ra\in M_n(R)$. This implies that $\mu$ is surjective, and therefore an isomorphism. Thus $K\otimes_R S_n(c,R)\simeq  S_n(c,K)$ as $K$-algebras.

Suppose that the $R$-algebras $S_n(c,R)$ and $S_m(d,R)$ are derived equivalent. Then there is a tilting complex $T$ for $S_n(c,R)$ such that $\End_{\Db{S_n(c,R)}}(T)\simeq S_m(d,R)$ as $R$-algebras. Since $K$ is a flat $R$-module, $\Tor^R_i(S_n(c,R),K)=0$ and $\Tor^R_i(S_m(d,R),K)=0$ for all $i\geq 1$. It then follows from \cite[Theorem 2.1]{JR2} that $K\otimes_R T$ is a tilting complex for $K\otimes_R S_n(c,R)$ with $\End_{\Db{K\otimes_R S_n(c,R)}}(K\otimes_R T)\simeq K\otimes_R S_m(d,R)$ as $K$-algebras.
Thus the $K$-algebras $S_n(c,K)$ and $S_m(d,K)$ are derived equivalent. By Theorem \ref{main1}, there holds $c\stackrel{D}{\sim}d$.
$\square$

\subsection{Relations among derived, Morita and stable equivalences: Proof of Corollary \ref{derp}}

Let $\Lambda$ be an Artin algebra and $M$ a generator-cogenerator for $\Lambda$-mod. Then the \emph{rigidity dimension} $\rd(M)$ of $M$ is defined by
$$\rd(M):=\operatorname {sup}\{n\in \mathbb{N}\mid {\Ext}^i_{\Lambda}(M,M)=0, \forall\ 1\leq i\leq n\}.$$
If no such $n$ exists, we define $\rd(M)=0$. The dominant dimension of the algebra $\Lambda$, denoted by $\dm(\Lambda)$, is the maximal $t\in \mathbb{N}$ (or $\infty$) such  that all the terms $I_0, I_1,\cdots, I_{t-1}$in a minimal injective resolution $$0\lra \Lambda \lra I_0 \lra I_1 \lra \cdots \lra I_t \lra \cdots$$ of $_{\Lambda}\Lambda$ are projective. By \cite[Lemma 3]{BJ}, $\dm(\End_\Lambda(M))=\rd(M)+2$.

The following lemma describes the dominant dimensions of principal centralizer matrix algebras.

\begin{Lem}\label{lem3.1}
$(1)$ $\dm(A_i)\in \{2,\infty\}$. Particularly, $\dm(S_n(c,R))\in \{2,\infty\}$.

$(2)$ $\dm(A_i)=\infty$ if and only if $A_i$ is a symmetric Nakayama algebra if and only if $P_c(f_i(x)^{n_i})$ is a singleton set. Thus $\dm(S_n(c,R))=\infty$ if and only if $S_n(c,R)$ is a symmetric Nakayama algebra if and only if $P_c(f_i(x)^{n_i})$ is a singleton set for all $i\in [l_c].$
\end{Lem}

{\it Proof.} If $\Lambda$ is an Artin algebra and $L\in \Lambda\modcat$, then it follows from the Auslander-Reiten formula $D\Ext^1_{\Lambda}(L,L)\simeq \ol{\Hom}_{\Lambda}(L,\tau L)$ that $\Ext^1_{\Lambda}(L,L)\neq 0$ if $\tau L\simeq L$, where $D$ is the usual duality of an Artin algebra and $\tau := D{\rm Tr}$ denotes the Auslander-Reiten translation.

Let $i\in [l_c]$. For the $U_i$-module $M_i$, $\tau (M_i)_{\mathscr{P}}\simeq (M_i)_{\mathscr{P}}$, and therefore $\rd(M_i)=\infty$ if $M_i$ is projective, and $0$, otherwise. Since $\dm(A_i)=\dm(\End_{U_i}(M_i))=\rd(M_i)+2,$ we deduce that $\dm(A_i)\in \{2,\infty\}$ and that $\dm(A_i)=\infty$ if and only if $M_i$ is projective if and only if $A_i$ is a symmetric Nakayama algebra if and only if $P_c(f_i(x)^{n_i})$ is a singleton set. Hence $\dm(S_n(c,R))\in \{2,\infty\}$. Moreover, $\dm(S_n(c,R))=\infty$ if and only if $S_n(c,R)$ is a symmetric Nakayama algebra if and only if $P_c(f_i(x)^{n_i})$ is a singleton set for all $i\in [l_c].$ $\square$

\medskip
Let $\Lambda$ be an Artin algebra. We denote by $\mathscr{P}(\Lambda)_{\mathscr{I}}$  the set of all isomorphism classes of projective $\Lambda$-modules without any nonzero injective summands.
In \cite[Proposition 1.5 and Theorems 1.7 and 2.6]{MV2}, Mart\'{i}nez-Villa proved the following.

\begin{Lem}\label{exa}{\rm \cite{MV2}} Let $F:\Lambda\stmc$ $\ra \bar{\Lambda}\stmc$ be a stable equivalence of Artin algebras $\Lambda$ and $\bar{\Lambda}$ both with neither nodes nor semisimple summands.

$(1)$ The functor $F$ provides a bijection $F':\mathscr{P}(\Lambda)_{\mathscr{I}}\ra \mathscr{P}(\bar{\Lambda})_{\mathscr{I}}$, which preserves simple projective modules in $\mathscr{P}(\Lambda)_{\mathscr{I}}$.

$(2)$ The functor $F$ induces a stable equivalence between the Frobenius parts of $\Lambda$ and $\bar{\Lambda}$.

$(3)$ Let $0\ra X\oplus Q_1\stackrel{f}\ra Y\oplus Q_2\oplus P\stackrel{g}\ra Z\ra 0$ be an exact sequence of $\Lambda$-modules without any split exact sequences as its direct summands, where $X,Y,Z\in \Lambda\modcat_{\mathscr{P}}$, $Q_1,Q_2\in \mathscr{P}(\Lambda)_{\mathscr{I}}$ and $P$ is a projective-injective $\Lambda$-module. Then there is a short exact sequence
$$0\lra F(X)\oplus F'(Q_1)\stackrel{f'}\lra F(Y)\oplus F'(Q_2)\oplus P'\stackrel{g'}\lra F(Z)\lra 0$$
in $\bar{\Lambda}\modcat$ such that $P'$ is projective-injective and that no split exact sequences are its direct summands.
\end{Lem}

We say that a stable equivalence $F:\Lambda\stmc{}\ra \bar{\Lambda}\stmc$ of Artin algebras
$\Lambda$ and $\bar{\Lambda}$ \emph{preserves non-semisimple blocks }if for non-projective indecomposable modules $M,N\in \Lambda\modcat_{\mathscr{P}}$, $F(M)$ and $F(N)$ lie in the same block of $\bar{\Lambda}$ if and only if $M$ and $N$ lie in the same block of $\Lambda$. In general, a stable equivalence may not preserve the numbers of non-semisimple blocks of algebras. This can be seen by \cite[Example 3.12]{xz3}, for instance.

In the following, we will show that stable equivalences between principal centralizer matrix algebras over a field do have this property.

As $S_n(c,R)\simeq\End_{R[c]}(R^n)\simeq \prod^{l_c}_{i=1}{\End}_{U_i}(M_i)=\prod^{l_c}_{i=1}A_i$, we consider the non-semisimple blocks of $S_n(c,R)$ and denote by $A$ the direct sum of its non-semisimple blocks. Now, we partition these blocks of $A$ in the following way such that

(1) $\mathcal{S}_{A,\geq 2}:=\{A_1, A_2, \cdots, A_{a_1}\}$ consists of the blocks $A_i$ with $n_i\ge 3$ and having at least 2 non-injective, indecomposable projective modules for $1\le i\le a_1$;

(2) $\mathcal{S}_{A,1}:=\{A_{a_1+1}, A_{a_1+2}, \cdots, A_{a_2}\}$ consists of the blocks $A_i$ with $n_i\ge 3$ and having only 1 non-injective, indecomposable projective module  for $a_1< i\le a_2$;

(3) $\mathcal{S}_{A,0}:=\{A_{a_2+1}, A_{a_2+2}, \cdots, A_{a_3}\}$ consists of the blocks $A_i$ with $n_i\ge 3$ and having  no non-injective, indecomposable projective modules for $a_2<i\le a_3$;

(4) $\mathcal{S}_{A}:=\{A_{a_3+1}, A_{a_3+2},\cdots, A_{a_4}\}$ consists of the blocks $A_i$ with $n_i=2$ for $a_3<i\le a_4$, where $0\le a_1\le a_2\le a_3\le a_4\le l_c.$

\smallskip
 Note that $n_i$ is the Loewy length of both $U_i$ and the center of $A_i$ by Lemma \ref{iso}. For an Artin algebra $\Lambda$, we have denoted by $\Lambda'$ the triangular matrix algebra obtained from $\Lambda$ by eliminating all nodes of $\Lambda$. Thus $\Lambda$ and $\Lambda'$ are stably equivalent, but the latter has no nodes.

Let $\widetilde{A}:=\prod^{a_3}_{i=1} A_i\; \times \prod_{a_3< i\le a_4} (A_i)'.$ Then we have the following result.

\begin{Lem}\label{replace} The algebra $\widetilde{A}$ has neither nodes nor semisimple direct summands. Moreover, there exists a stable equivalence
$F_A: A\stmc{}\ra \widetilde{A}\stmc$ such that $F_A$ preserves non-semisimple blocks of algebras and that the restriction of $F_A$ to the block $A_i$ is induced by the identity functor for all $i\in [a_3]$.
\end{Lem}

{\it Proof.} Recall that $U_i= R[x]/(f_i(x)^{n_i})$ and $A_i = \End_{U_i}(M_i)$. For $i\in [a_3],$ it follows by Lemma \ref{node} that the block $A_i$ has neither nodes nor projective simple modules. For $a_3<j\le a_4$, the triangular matrix algebra $(A_j)'$ is stably equivalent to $A_j$ and has no nodes. Clearly, $\widetilde{A}$ does not have any semisimple direct summands. Now it is easy to get a desired stable equivalence $F_A$ between $A$ and $\widetilde{A}$. $\square$

\medskip
Let $S_m(d,R)$ be another centralizer matrix algebra and $B$ be the sum of its non-semisimple blocks. Similarly, we have a partition of blocks for $B$. This is given by the natural numbers $0\le b_1\le b_2\le b_3\le b_4\le l_d$, namely the partition $\{B_1, \cdots, B_{b_1}\}\cup \{ B_{b_1+1},\cdots, B_{b_2}\} \cup \{B_{b_2+1},\cdots,B_{b_3}\}\cup \{B_{b_3+1},\cdots, B_{b_4}\}$ has the corresponding properties as the blocks of $A$.

By Lemma \ref{replace},
$\widetilde{B}:=\prod^{b_3}_{j=1} B_j \; \times \prod_{b_3< j\le b_4} (B_j)'$
has neither nodes nor semisimple direct summands, and there is a stable equivalence $F_B: B\stmc{}\ra \widetilde{B}\stmc$ preserving non-semsimple blocks such that the restriction of $F_B$ to $B_j$ is induced by the identity functor for $j\in [b_3]$.

Now, we assume that there is a stable equivalence $F$ between $S_n(c,R)$ and $S_m(d,R)$. Then $F$ restricts to a stable equivalence between $A$ and $B$. Thus $H:=F_B\circ F\circ F_A^{-1}:\widetilde{A}\stmc$ $\to$ $\widetilde{B}\stmc$ is a stable equivalence. Let $J$ be a quasi-inverse of $H$. As defined in Lemma \ref{exa}(1), $H': \mathscr{P}(\widetilde{A})_{\mathscr{I}}\ra \mathscr{P}(\widetilde{B})_{\mathscr{I}}$ and $J': \mathscr{P}(\widetilde{B})_{\mathscr{I}}\ra \mathscr{P}(\widetilde{A})_{\mathscr{I}}$ are the bijections induced by $H$ and $J$, respectively.

\begin{Lem}\label{bij}
The correspondence $H'$ induces a bijection between $\mathcal{S}_{A,\geq 2}\cup \mathcal{S}_{A,1}$ and $\mathcal{S}_{B,\geq 2}\cup \mathcal{S}_{B,1}$ such that the corresponding blocks have the same number of non-injective, indecomposable projective modules. In particular, $a_l=b_l$ for $1\le l\le 2$.
\end{Lem}

{\it Proof.} By the proof in Lemma \ref{node}, for $i\in [a_2]$, the block $A_i$ of $\widetilde{A}$ has a unique projective-injective indecomposable module, say $P_i$. We have to discuss the following $2$ cases.

(i) $H$ induces a bijection from  $\mathcal{S}_{A,\geq 2}$ to $\mathcal{S}_{B,\geq 2}$, such that the corresponding blocks have the same number of non-injective, indecomposable projective modules.

In fact, let $A_i\in \mathcal{S}_{A,\geq 2}$, that is, $A_i$ has at least $2$ non-injective, indecomposable projective modules $P_{i1}$ and $P_{i2}$ which are non-simple by Lemma \ref{node}. Then there exist $2$ indecomposable direct summands $M_{i1}$ and $M_{i2}$ of the $U_i$-module $M_i$ such that $P_{ir}\simeq \Hom_{U_i}(M_i,M_{ir})$ as $A_i$-modules for $1\le r\le 2$. As $U_i$ is a local Nakayama algebra, we may assume that $M_{i1}$ is isomorphic to a proper submodule of $M_{i2}$. It then follows from the left exactness of the Hom-functor $\Hom_{U_i}(M_i,-)$ that $P_{i1}$ is isomorphic to a submodule of $P_{i2}$. Hence there is
an exact sequence of $A_i$-modules
$$0\lra P_{i1}\lra P_{i2}\lra P_{i2}/P_{i1}\lra 0$$
with $P_{i2}/P_{i1}$ indecomposable. By Lemma \ref{exa}(3), there is an exact sequence of $\widetilde{B}$-modules
$$(*)\quad 0\lra H'(P_{i1})\lra H'(P_{i2})\oplus P'_1\lra H(P_{i2}/P_{i1})\lra 0,$$
with $P'_1$ being projective-injective, such that this sequence does not contain split exact sequences as its direct summands.
Then $\Hom_{\widetilde{B}}(H'(P_{i1}), H'(P_{i2}))\ne 0$, and therefore $H'(P_{i1})$ and $H'(P_{i2})$ lie in the same block of $\widetilde{B}$. Otherwise,  the sequence $(*)$ would contain a split sequence $0\to 0 \to H'(P_{i2})\stackrel{1}{\to} H'(P_{i2})\to 0$ as its summand. Let $O$ be the block of $\widetilde{B}$ to which $H'(P_{i1})$ and $H'(P_{i2})$ belong. Note that a block $(B_j)'$ of $\widetilde{B}$, with $b_3<j\le b_4$, has at most $2$ non-injective, indecomposable projective modules, one of which is a simple module by Remark \ref{rmk2.11}. It then follows from Lemma \ref{exa}(1) that $O\not\simeq (B_j)'$ as algebras for $b_3<j\le b_4$. Hence $O$ is a block of the form $B_j$ for $j\in [b_1]$ and $H'$ sends all the non-injective, indecomposable projective $A_i$-modules to the ones belonging to the block $B_j.$  Similarly, the non-injective, indecomposable projective $B_j$-modules are mapped by $J'$ into modules belonging to the block $A_i.$ Thus $H'$ restricts to a bijection between the set of non-injective, indecomposable projective $A_i$-modules and the one of non-injective, indecomposable projective $B_j$-modules. This implies that $H$ induces a bijection from the set of the blocks in $\mathcal{S}_{A,\geq 2}$ to the set of blocks in $\mathcal{S}_{B,\geq 2}$, such that the corresponding blocks have the same number of non-injective, indecomposable projective modules. Clearly, $a_1=b_1$.

(ii) $H$ induces a bijection between $\mathcal{S}_{A,1}$ and $\mathcal{S}_{B,1}$.

Actually, let $A_i\in \mathcal{S}_{A,1}$, that is, $A_i$ has only $1$ non-injective, indecomposable projective module, say $Q_i$. If $H'(Q_i)$ lies in some block $(B_j)'$ for $b_3<j\le b_4$, then, by Remark \ref{rmk2.11}, $H'(Q_i)$ has a simple projective submodule, say $\widetilde{P_j}$. With a similar argument as in (i), we deduce that $Q_i$ and the simple projective module $J'(\widetilde{P_j})$ lie in the same block $A_i$. Note that $Q_i$ is not simple and $Q_i\not\simeq J'(\widetilde{P_j})$. This implies that the block $A_i$ contains $2$ non-injective, indecomposable projective modules, a contradiction. Thus it follows from (i) that $H'(Q_i)$ belongs to a block $B_j\in \mathcal{S}_{B,1}$. So $H$ induces a bijection from the set of blocks in $\mathcal{S}_{A,1}$ to  the set of blocks in $\mathcal{S}_{B,1}$, and therefore $a_2-a_1=b_2-b_1$ and $a_2=b_2$. $\square$

\begin{Lem}\label{p-i} Let $\widetilde{A_1}$ and $\widetilde{B_1}$ be the sum of blocks in $\mathcal{S}_{A,\geq 2}\cup \mathcal{S}_{A,1}$ and $\mathcal{S}_{B,\geq 2}\cup \mathcal{S}_{B,1}$, respectively.
Then the functor $H$ restricts to a stable equivalence between $\widetilde{A_1}$ and $\widetilde{B_1}$ preserving non-semisimple blocks.
\end{Lem}

{\it Proof.} By Lemma \ref{bij}, we assume $H'(A_i)=B_i$ for $i\in [a_2]$. Let $i\in [a_2]$ and $A_i$ be a block in $\widetilde{A_1}$. Suppose that $M$ is a non-projective, indecomposable $A_i$-module. Then $H(M)$ is indecomposable. Further, we show that $H(M)$ lies in the block $B_i$ of $\widetilde{B_1}$.

(a) By the proof of Lemma \ref{node}, any non-injective, indecomposable projective $A_i$-module is isomorphic to a submodule of the unique projective-injective $A_i$-module $P_i$. This implies that $\rad(P_i)$ is indecomposable and each non-injective, indecomposable projective $A_i$-module is isomorphic to a submodule of $\rad(P_i)$. If $\rad(P_i)$ is projective, then $H'(\rad(P_i))$ lies in the block $B_i$. If $\rad(P_i)$ is not projective, then there exists a non-injective, indecomposable projective $A_i$-module $P$ and an exact sequence
$$\quad 0\lra P\stackrel{\iota}\lra \rad(P_i)\stackrel{\eta}\lra \rad(P_i)/P\lra 0$$
without split direct summands.
Applying Lemma \ref{exa}(3) to this sequence, we get an exact sequence of $\widetilde{B}$-modules
$$0\lra H'(P)\stackrel{\iota'}\lra H(\rad(P_i))\oplus P'_2\stackrel{\eta'}\lra H(\rad(P_i)/P)\lra 0,$$
where $P'_2$ is projective-injective. We show that $H(\rad(P_i))$ lies in the block $B_i.$ Suppose contrarily  that $H(\rad(P_i))$ does not belong to the block $B_i.$ Then $\Hom_{\widetilde{B}}(H'(P), H(\rad(P_i)))=0$, and therefore $\Img(\iota')\subset P'_2.$ Thus $H(\rad(P_i)/P)\simeq (H(\rad(P_i))\oplus P'_2)/\Img(\iota')\simeq P'_2/\Img(\iota')\oplus H(\rad(P_i))$. This yields that $\rad(P_i)$ is isomorphic to a direct summand of $\rad(P_i)/P$, a contradiction.

(b) Let $P(M)$ be the projective cover of $M$. Then there is an exact sequence of $A_i$-modules
$$ 0\lra \Omega(M)\stackrel{f}\lra P(M)\stackrel{g}\lra M\lra 0,$$
which has no split direct summands.
Write $\Omega(M)=L_1\oplus L_2$ for $L_1\in {A_i\modcat}_{\mathscr{P}}$ and $L_2\in \mathscr{P}(A_i)_{\mathscr{I}}$, and $P(M)=Q\oplus (P^{\oplus t}_i)$ with $Q\in \mathscr{P}(A_i)_{\mathscr{I}}$ and $t\in \mathbb{N}.$ By Lemma \ref{exa}(3), we get an exact sequence of $\widetilde{B}$-modules
$$(\ddag)\quad 0\lra H(L_1)\oplus H'(L_2)\stackrel{f'}\lra H'(Q)\oplus P'_3\stackrel{g'}\lra H(M)\lra 0,$$
such that $P'_3$ is projective-injective and ($\ddag$) has no split direct summands.

(1)Assume $Q\neq 0$. Then $H(M)$ has to lie in $B_i$. Otherwise, we would have $\Hom_{\widetilde{B}}(H'(Q),H(M))=0$, and therefore $ \Img(f')=H'(Q)\oplus (\Img(f')\cap P'_3).$ Since $H'(Q)$ is projective, there exists a homomorphism $h':H'(Q)\ra H(L_1)\oplus H'(L_2)$ such that $h'f'$ is the identity on $H'(Q)$, and therefore
$$0\lra h'(H'(Q))\stackrel{f'}\lra H'(Q)\lra 0\lra 0$$is a split direct summand of ($\ddag$), a contradiction. Hence $H(M)$ belongs to the block $B_i$.

(2) Assume $Q=0$. Then $H'(Q)=0.$ Suppose that $P'_3=\bigoplus^l_{k=1} C_k$, where $l\ge 1$ and $C_k$ is a sum of indecomposable direct summands of $P'_3$ belonging to the same block. We shall show that all indecomposable direct summands of $P'_3$ lie in the same block of $\widetilde{B}$, that is, $l=1$.
Now, suppose contrarily that $l\geq 2$. By $(\ddag)$, $\Img(f')=\bigoplus^l_{k=1} D_k$, where $D_k$ is a submodules of $C_k$ for $k\in [l]$. Since $H(L_1)\oplus H'(L_2)$ contains no projective-injective direct summands, each $D_k$ is a proper submodules of $C_k$. Therefore $H(M)\simeq P'_3/\Img(f')\simeq \bigoplus^l_{k=1} C_k/D_k$ is decomposable, a contradiction. Thus $l=1$ and all indecomposable direct summands of $P'_3$ lie in the same block of $\widetilde{B}$. This also imply that $P'_3$ belongs to the same block.

Suppose that $\Omega(M)$ has a direct summand isomorphic to $\rad(P_i).$ Then $f'$ in $(\ddag)$ restricts to an injective homomorphism from  $H(\rad(P_i))$ if $\rad(P_i)$ is not projective (or from $H'(\rad(P_i))$ if $\rad(P_i)$ is projective) to $P'_3$. In particular, $H(\rad(P_i))$ (or $H'(\rad(P_i))$) and $P'_3$ lie in the same block $B_i$ of $\widetilde{B}$. Also, $H'(Q)=0$ implies $\Hom_{\widetilde{B}}(P'_3,H(M))\neq 0$ in $(\ddag)$. This implies that $H(M)$ and $P'_3$ lie in the same block $B_i$.

Suppose that $\Omega(M)$ has no direct summands isomorphic to $\rad(P_i).$ Recall that $P(M)=P_i^{\oplus t}$ under the assumption $Q=0$. We consider the exact sequence of $A_i$-modules
$$(\sharp)\quad 0\lra \Omega(M)\stackrel{f}\lra (\rad(P_i))^{\oplus t}\stackrel{g}\lra \rad(M)\lra 0.$$
Deleting the split direct summands of ($\sharp$), we obtain an exact sequence of $A_i$-modules
$$0\lra \Omega(M)\stackrel{f_0}\lra (\rad(P_i))^{\oplus r}\stackrel{g_0}\lra X\lra 0$$
for some $r\le t$ and a submodule $X$ of $\rad(M)$.
Thanks to Lemma \ref{exa}(3), there is an indecomposable direct summand $L$ of $\Omega(M)$ such that $H(L)$ (or $H'(L))$ lies in the block $B_i$. Otherwise $H(X)$ would contain a direct summand isomorphic to $H((\rad(P_i))^{\oplus r})$ (or $H'((\rad(P_i))^{\oplus r})$) (see the argument in (a)), a contradiction to that $X$ contains no direct summands isomorphic to $(\rad(P_i))^{\oplus r}).$ Now, we see from ($\ddag$) that the modules $H(M)$, $P'_3$ and $H(L)$ lie in the same block $B_i$ of $\widetilde{B}$ or  the modules $H(M)$, $P'_3$ and $H'(L)$ lie in the same block $B_i$ of $\widetilde{B}$.

Thus we have proved that $H(M)$ lies in the block $B_i$ of $\widetilde{B_1}$.
Similarly, for $i\in [a_2]$ and a non-projective, indecomposable $B_i$-module $N$, we see that $J(N)$ belongs to the block $A_i$ of $\widetilde{A_1}$. Hence $H$ induces a stable equivalence between $\widetilde{A_1}$ and  $\widetilde{B_1}$, which preserves non-semisimple blocks. $\square$

\medskip
For an Artin algebra $\Lambda$, let $\Gamma_{\Lambda}$ denote the Auslander-Reiten quiver of $\Lambda$ and $\Gamma^s_{\Lambda}$ the stable Auslander-Reiten quiver of $\Gamma_{\Lambda}$ obtained by removing all projective vertices from $\Gamma_{\Lambda}$. For a local, symmetric Nakayama algebra $\Lambda_0:= R[x]/(f(x)^n)$, $\Gamma^s_{\Lambda_0}$ is a connected quiver such that there are two arrows between any two vertices if they are connected by an irreducible map.

\begin{Lem}\label{no-p-i}Let $\widetilde{A_2}$ and $\widetilde{B_2}$ be the sum of blocks in $\mathcal{S}_{A,0}\cup \mathcal{S}_A$ and $\mathcal{S}_{B,0}\cup  \mathcal{S}_B$, respectively.
The functor $H$ restricts to a stable equivalence between $\widetilde{A_2}$ and $\widetilde{B_2}$, which preserves non-semisimple blocks.
\end{Lem}

{\it Proof.} Given a non-projective, indecomposable $\widetilde{A_2}$-module $Y$, the module $H(Y)$ belongs to a block of $\widetilde{B_2}$ by Lemma \ref{p-i}. Let $A_i\in \mathcal{S}_{A,0}$ be a block, that is, all indecomposable projective $A_i$-modules are injective. Then the $U_i$-module $M_i$ is projective and $A_i= \End_{U_i}(M_i)$ is a symmetric Nakayama algebra which is Morita equivalent to $U_i=R[x]/(f_i(x)^{n_i})$. We show that $H$ restricts to a stable equivalence between $A_i$ and a block in $\mathcal{S}_{B,0}.$ Let $\{K_l\mid l\in [n_i-1]\}$ be the set of all (up to isomorphism) non-projective, indecomposable $A_i$-modules. Since $\Gamma^s_{A_i}$ is a connected quiver with two arrows between any two connected vertices, it follows from \cite[Lemma 1.2(d), p. 336]{ARS} that all modules $H(K_l), l\in [n_i-1]$, lie in the same Auslander-Reiten component of a block $W$ in $\widetilde{B_2}$. For a block $B_k$ in $\mathcal{S}_B$, we consider the node-eliminated block $(B_k)'$. Since $B_k$ and $(B_k)'$ are stably equivalent,  the quiver $\Gamma^s_{B_k}$ and $\Gamma^s_{(B_k)'}$ are isomorphic as translation quivers by \cite[Corollary 1.10, p.342]{ARS}. According to Remark \ref{rmk2.11}, for $b_3<k\le b_4$, $\Gamma^s_{B_k}$ either contains only $1$ vertex or is of the form $\bullet\ra\bullet\ra\bullet$. This implies that $W$ can not be a block $(B_k)'$ in $\widetilde{B_2}$ with $b_3<k\le b_4.$ Thus $W$ is a block in $\mathcal{S}_{B,0}$, say $B_j$. Clearly, $B_j$ is a symmetric Nakayama algebra and all $H(K_l), l\in [n_i-1]$, are precisely the non-projective, indecomposable $B_j$-modules (up to isomorphism). Similarly, for each block $B_j\in \mathcal{S}_{B,0}$, there is a unique block $A_i$ in $\mathcal{S}_{A,0}$ such that $J$ restricts to a stable equivalence between $B_j$ and $A_i$. In this way, $H$ induces not only a one-to-one correspondence but also a stable equivalence between the blocks in $\mathcal{S}_{A,0}$ and $\mathcal{S}_{B,0}.$

Let $A_i\in \mathcal{S}_A$ be a block, that is, $a_3<i\le a_4.$ The quiver $\Gamma^s_{(A_i)'}$ for $a_3<i\le a_4$ (respectively, $\Gamma^s_{(B_j)'}$ for $b_3<j\le b_4$) is connected with either $1$ or $3$ vertices. Clearly, $H$ restricts to a stable equivalence between $(A_i)'$ and some block $(B_j)'$ with $b_3<j\le b_4$. Hence $H$ induces not only a one-to-one correspondence but also a stable equivalence between the the blocks in $\mathcal{S}_A$ and $\mathcal{S}_B.$ $\square$

\begin{Lem}\label{sta}
Let $c\in M_n(R)$ and $d\in  M_m(R)$. Suppose that there is a stable equivalence $F$ between $S_n(c,R)$ and $S_m(d,R)$. Then $F$ preserves non-semisimple blocks. Moreover, if $A_i$ and $B_j$ are stably equivalent, then $n_i=m_j$.
\end{Lem}

{\it Proof.} We keep all notions introduced previously. Suppose that there is a stable equivalence $F$ between $S_n(c,R)$ and $S_m(d,R)$.
Recall that $A$ and $B$ are the sum of non-semisimple blocks of $S_n(c,R)$ and $S_m(d,R)$, respectively, and $F_A: A\stmc{}\ra \widetilde{A}\stmc$ and $F_B: B\stmc{}\ra \widetilde{B}\stmc$ are stable equivalences.

By Lemmas \ref{p-i} and \ref{no-p-i}, the stable equivalence $H=F_B\circ F\circ F_A^{-1}:\widetilde{A}\stmc{}\ra \widetilde{B}\stmc$ preserves non-semisimple blocks. Since $F_A$ and $F_B$ preserve non-semisimple blocks, we infer  that $F$ preserves non-semisimple blocks.

Suppose that the blocks $A_i$ and $B_j$ are stably equivalent. Then  $A_i$ is semisimple if and only if $B_j$ is semisimple. In this case, $n_i=m_j=1$. So we may assume that $n_i\geq 2$ and $m_j\geq 2$. We first show that $n_i=2$ if and only if $m_j=2$. Suppose contrarily that either $n_i=2$ and $m_j\geq 3$ or $n_i\geq 3$ and $m_j=2$. We only deal with the situation $n_i=2$ and $m_j\geq 3$. The other case can be done similarly. It follows from Lemma \ref{node} that  $A_i$ has nodes but $B_j$ does not have nodes. We replace $A_i$ by a stably equivalent algebra $(A_i)'$ without nodes. By Lemma \ref{exa}(2), the Frobenius parts of $(A_i)'$ and $B_j$ are stably equivalent, while the Frobenius parts of $(A_i)'$ is zero by Remark \ref{rmk2.11} and the Frobenius parts of $B_j$ is $V_j$ by the proof of Theorem \ref{main1}(4). Thus $V_j\stmc$ is zero. This is a contradiction. Thus $n_i=2$ if and only if $m_j=2$.
Now assume  $n_i\geq 3$ and $m_j\geq 3.$ With a similar argument, we deduce that the Frobenius parts of $A_i$ and $B_j$ are stably equivalent, that is, $U_i$ and $V_j$ are stably equivalent. In particular, $U_i$ and $V_j$ have the same number of non-projective, indecomposable modules, that is, $n_i-1=m_j-1.$ Hence $n_i=m_j$.$\square$

\medskip
{\bf Proof of Corollary \ref{derp}}. We keep all notion introduced previously. Let $c\in M_n(R)$ and $d\in M_m(R)$.

(1) Assume that $c$ and $d$ are permutation matrices and that $S_n(c,R)$ and $S_m(d,R)$ are derived equivalent. Then $S_n(c,R)$ and $S_m(d,R)$ have the same number of blocks, that is, $l_c=l_d$. So we may assume that $A_i$ and $B_i$ are derived equivalent for $i\in [l_c].$ By Lemma \ref{3.1}, $U_i\simeq V_i$ and $n_i=m_i$ for $i\in [l_c].$ By Theorem \ref{main1}(1), it suffices to show that $P_c(f_i(x)^{n_i})=P_d(g_i(x)^{m_i})$ for $i\in [l_c].$ By Lemma \ref{per}, the integers in $P_c(f_i(x)^{n_i})$ are $p$-powers for $i\in [l_c]$. Similarly, the integers in $P_c(g_i(x)^{m_i})$ are $p$-powers for $i\in [l_c]$. Let $t_i:=|P_c(f_i(x)^{n_i})|=|P_d(g_i(x)^{n_i})|$ for $i\in [l_c]$. If $t_i=1$, then $P_c(f_i(x)^{n_i})=\{n_i\}=\{m_i\}=P_d(g_i(x)^{m_i})$. For instance, if $p=0$, then $t_i=1$. Now, we may assume that $t_i\geq 2$ and $p>0$. Let $P_c(f_i(x)^{n_i}):=\{p^{u_1},\cdots,p^{u_{t_i}}\}$ with $u_1> \cdots > u_{t_i}$ and $P_d(g_i(x)^{m_i}):=\{p^{v_1},\cdots,p^{v_{t_i}}\}$ with $v_1> \cdots > v_{t_i}$. By Theorem \ref{main1}(2), we get $\{p^{u_1}-p^{u_{2}},\cdots, p^{u_{t_i-1}}-p^{u_{t_i}},p^{u_{t_i}}\}= \{p^{v_1}-p^{v_2},\cdots, p^{v_{t_i-1}}-p^{v_{t_i}},p^{v_{t_i}}\}.$
For positive integers $a>b$ and $s>t$, the number $p^a-p^b$ is a $p$-power if and only if $p=2$ and $a=b+1$; and the equality $p^a-p^b=p^s-p^t$ holds if and only if $a=s$ and $b=t$. By considering the cases $p=2$  and $p\geq 3$ separately, we get $u_k=v_k$ for all $k\in [t_i].$  Thus $P_c(f_i(x)^{n_i})=P_d(g_i(x)^{m_i})$ for $i\in [l_c].$ This implies that $A$ and $B$ are Morita equivalent by Theorem \ref{main1}(1).

(2) Let $R$ be a  perfect field. Suppose that $S_n(c,R)$ and $S_m(d,R)$ are representation-finite and stably equivalent. Let $F$ be a stable equivalence between $S_n(c,R)$ and $S_m(d,R).$ Then it follows from Lemma \ref{rep-f} that
$$P_c(f_i(x)^{n_i})\subseteq \{1,\max\{n_i,3\}-1,\max\{n_i,3\}\} \mbox{ and  } P_d(g_j(x)^{m_j})\subseteq \{1,\max\{m_j,3\}-1,\max\{m_j,3\}\}$$for all $f_i(x)^{n_i}\in \mathcal{M}_c$ and $g_j(x)^{m_j}\in \mathcal{M}_d.$ Reordering the blocks of $S_n(c,R)$ and $S_m(d,R)$, we may assume that there are natural numbers $a\le a'$ and $b\le b'$ such that $n_i\geq 3$ (respectively, $m_j\geq 3$) if and only if $i\in [a]$ (respectively, $j\in [b]$), and that $n_i=2$ (respectively, $m_j=2$) if and only if $a<i\le a'$ (respectively, $b<j\le b'$). Let $A$ and $B$ be the sums of non-semisimple blocks in $S_n(c,R)$ and $S_m(d,R)$, respectively. By Lemma \ref{sta}, $F$ preserves non-semisimple blocks. Moreover, if $A_i$ and $B_j$ are stably equivalent, then $n_i=m_j$. Thus $a=b$ and $a'=b'.$ Therefore  we may assume that $F$ restricts to a stable equivalence $F_i$ between $A_i$ and $B_i$ for $i\in [a']$. In particular, $n_i=m_i$ for $i\in [a'].$

Let $i\in [a].$ Then $n_i=m_i\ge3$. By Lemma \ref{node}(2), $A_i$ and $B_i$ have no nodes. By Lemma \ref{exa}(2), $F_i$ induces a stable equivalence between the Frobenius parts of $A_i$ and $B_i$, that is, $U_i$ and $V_i$ are stably equivalent. Let $K_i:=U_i/\rad(U_i)$. Since $f_i(x)$ is separable by our assumption on the ground field, it follows from Corollary \ref{St-i} that $U_i\simeq K_i[x]/(x^{n_i})\simeq V_i$ as algebras.

According to Lemma \ref{exa}(1), $A_i$ and $B_i$ have the same number of non-injective, indecomposable projective modules, and hence the same number of indecomposable projective modules. It then follows from $\mathcal{B}(M_i)\simeq \bigoplus_{r\in {P_c(f_i(x)^{n_i})}} R[x]/(f_i(x)^r)$ as $U_i$-modules and $\mathcal{B}(N_j)\simeq \bigoplus_{s\in {P_d(g_j(x)^{m_j})}} R[x]/(g_j(x)^s)$ as $V_j$-modules that $$|P_c(f_i(x)^{n_i})|=|P_d(g_i(x)^{n_i})|.$$Thus, due to the inclusions $\{n_i\}\subseteq P_c(f_i(x)^{n_i})\subseteq \{1,n_i-1,n_i\}$ and $\{n_i\}\subseteq P_d(g_i(x)^{n_i})\subseteq \{1,n_i-1,n_i\}$, we obtain $P_c(f_i(x)^{n_i})=P_d(g_i(x)^{n_i})$ or $P_c(f_i(x)^{n_i})=\mathcal{J}_{P_d(g_i(x)^{n_i})}.$ Hence $A_i$ and $B_i$ are almost $\nu$-stable derived equivalent by Lemma \ref{alm} if we identify $U_i$ with $V_i$, and therefore $A_i$ and $B_i$ are stably equivalent of Morita type by \cite[Theorem 1.1]{hx1}.

Suppose $a<i\le a'.$ Then $n_i=m_i=2$. As the argument in the foregoing case $i\in [a]$, we can show that  $U_i\simeq V_i$ as algebras and that $A_i$ and $B_i$ are Morita equivalent.
Thus $A$ and $B$ are stably equivalent of Morita type. Since $R$ is perfect, the semisimple blocks of $S_n(c,R)$ and $S_m(d,R)$ are separable $R$-algebras. Hence $S_n(c,R)$ and $S_m(d,R)$ are stably equivalent of Morita type.

(3) Suppose that $S_n(c,R)$ and $S_m(d,R)$ are derived equivalent. By Lemma \ref{lem3.1}(1), $\dm(S_n(c,R)) \in \{2,\infty\}$. Thus, to prove that $S_n(c,R)$ and $S_m(d,R)$ have the same dominant dimension, we show that $\dm(S_n(c,R)) = \infty$ if and only if $\dm(S_m(d,R))=\infty.$ However, this follows from Theorem \ref{main1}(2) and Lemma \ref{lem3.1}(2) immediately. Thus $S_n(c,R)$ and $S_m(d,R)$ have the same dominant dimension.

Now, suppose that $S_n(c,R)$ and $S_m(d,R)$ are stably equivalent. Assume $\dm(S_n(c,R))=\infty.$ Then $S_n(c,R)$ is a symmetric Nakayama algebra by Lemma \ref{lem3.1}(2). By \cite[Corollary 1.2]{IR}, every non-simple projective $S_m(d,R)$-module is injective, while the indecomposable projective modules are of the form $\Hom_{V_j}(N_j,R[x]/(g_j(x)^s))$ for $s\in P_d(g_j(x)^{m_j})$ and $j\in [l_d]$. Thus $P_d(g_j(x)^{m_j})$ is a singleton set for all $g_j(x)^{m_j}\in \mathcal{M}_d$, and therefore $\dm(S_m(d,R))=\infty$ by Lemma \ref{lem3.1}(2).
$\square$

\medskip
For representation-finite, self-injective algebras over an algebraically closed field, Asashiba proved in \cite{As} that stable equivalences lift to stable equivalences of Morita type. His proof uses classification of representation-finite, self-injective algebras under derived equivalences. In general, principal centralizer matrix algebras do not have to be self-injective. As shown in the above, our proof uses a completely different strategy.

\begin{Koro}\label{one-more}
Let $R$ be a noetherian domain of characteristic $p>0$ and $\sigma\in \Sigma_n$ be of cycle type $\lambda:=(\lambda_1, \cdots, \lambda_s)$, and let $\sigma^+$ be a permutation in $\Sigma_{n+1}$ of cycle type $\lambda^+:=(\lambda_1, \cdots, \lambda_s, 1)$. Then the following are equivalent

$(a)$ $S_n(c_\sigma,R)$ and $S_{n+1}(c_{\sigma^+},R)$ are derived equivalent.

$(b)$ $S_n(c_\sigma,R)$ and $S_{n+1}(c_{\sigma^+},R)$ are Morita equivalent.

$(c)$ There is an $i\in [s]$ such that $p\nmid \lambda_i$.
\end{Koro}

{\it Proof.} Let $K$ be the fraction field of $R$ and $F_p$ be the prime field of $K$. Since $c_{\sigma^+}$ is just the diagonal block-matrix diag$(c_{\sigma},1)$, we have $\mathcal{E}_{c_{\sigma^+}}=\mathcal{E}_{c_\sigma}\cup \{x-1\}$ when $c_\sigma$ and $c_{\sigma^+}$ are viewed as matrices over either $K$ or $F_p$.
Note that all $\lambda_i$ are exactly the orbit lengths of $\langle\sigma\rangle$ on $[n]$.

$(a)\Rightarrow (c)$ Suppose $S_n(c_\sigma,R)$ and $S_{n+1}(c_{\sigma^+},R)$ are derived equivalent. Then it follows from Remark \ref{rmk3.3} that $S_n(c_\sigma,K)$ and $S_{n+1}(c_{\sigma^+},K)$ are derived equivalent, and hence Morita equivalent by Corollary \ref{derp}. Further, by Lemma \ref{per}, $p\nmid \lambda_i$ for some number $i$.

$(c)\Rightarrow (b)$ Assume (c). Then it follows from Lemma \ref{per} that $x-1\in \mathcal{E}_{c_\sigma}.$ By Theorem \ref{main1}, $S_n(c_\sigma,F_p)$ and $S_{n+1}(c_{\sigma^+},F_p)$ are Morita equivalent. With an argument similar to the one in Remark \ref{rmk3.3}, we obtain the isomorphisms of $R$-algebras
$$R\otimes_{F_p}S_n(c_\sigma,F_p)\simeq S_n(c_\sigma,R)\; \mbox{ and } \; R\otimes_{F_p}S_m(c_{\sigma^+},F_p)\simeq S_m(c_{\sigma^+},R).$$ Hence $S_n(c_\sigma,R)$ and $S_{n+1}(c_{\sigma^+},R)$ are Morita equivalent.
$\square$

\medskip
Finally, we consider the case of nilpotent matrices. Let $c\in M_n(R)$ be a nilpotent matrix. Then the Jordan canonical form $c_0$ of $c$ is unique up to the ordering of its Jordan blocks. Further, $c_0$ has a Jordan block of size $t$ if and only if ${\rm rank}(c^{t+1})+{\rm rank}(c^{t-1})-2{\rm rank}(c^{t})>0$. We set $I_c:=\{t\ge 1\mid c_0~\mbox{ has a Jordan block of size t}\}.$ Note that $\mathcal{M}_c$ consists of only one polynomial of the form $x^r$ with $r$ being the maximal number in $I_c$. Thus $I_c = P_c(x^r).$

\begin{Koro}
Let $c\in M_n(R)$ be a nilpotent matrix and $d\in M_m(R)$. Then $S_n(c,R)$ and $S_m(d,R)$ are derived equivalent if and only if $d=\lambda I_m+b$ with $\lambda\in R$ and $b$ being a nilpotent matrix such that $\mathcal{H}_{I_b}=\mathcal{H}_{I_c}.$
\end{Koro}

{\it Proof.} Suppose $d=\lambda I_m +b$ with $\lambda \in R$ and $b\in M_m(R)$ a nilpotent matrix, such that $\mathcal{H}_{I_b}=\mathcal{H}_{I_c}.$ Then $S_m(d,R)=S_m(b,R).$ Let $x^s$ be the unique polynomial in  $\mathcal{M}_b$. Clearly, the assumption $\mathcal{H}_{I_c}=\mathcal{H}_{I_b}$ implies $\mathcal{H}_{P_c(x^r)}=\mathcal{H}_{P_b(x^s)}.$ It then follows from Theorem \ref{main1}(2) that $S_n(c,R)$ and $S_m(b,R)$ are derived equivalent.

Conversely, suppose that $S_n(c,R)$ and $S_m(d,R)$ are derived equivalent. Then it follows from Theorem \ref{main1}(2) that $\mathcal{M}_d$ consists of only one polynomial, say $f(x)^s$ with an irreducible polynomial $f(x)\in R[x]$ and $s\in \mathbb{N}$, and that $R[x]/(x^r)\simeq R[x]/(f(x)^s)$ as algebras. Thus $r=s$ and $f(x)=x-\lambda$ for some $\lambda\in R$. Set $b:=\lambda I_m- d$. Then $m_{b}(x)=x^s$, that is, $b$ is a nilpotent matrix. Clearly, $P_d(f(x)^s)=P_b(x^s).$ Therefore $\mathcal{H}_{P_c(x^r)}=\mathcal{H}_{P_d(f(x)^s)}=\mathcal{H}_{P_b(x^s)}.$ Hence $\mathcal{H}_{I_c}=\mathcal{H}_{I_b}.$ $\square$

\subsection{Restrictions of derived equivalences\label{sect4}}
In this section we investigate the relation between a derived equivalence of centralizer algebras of permutation matrices and the one of their $p$-parts.

Let $\sigma=\sigma_1\cdots\sigma_s\in \Sigma_n$ be the product of disjoint cycle-permutations $\sigma_i$ of cycle type $\lambda=(\lambda_1,\cdots, \lambda_s)$ with $\lambda_i\ge 1$ for $1\le i\le s$. For a prime number $p>0$,  a cycle $\sigma_i$ is said to be \emph{$p$-regular} if $p\nmid \lambda_i$, and \emph{$p$-singular} if $p\mid \lambda_i$. If $p=0$, all cycles are $p$-regular. Let $r(\sigma)$ (respectively, $s(\sigma)$) be the product of the $p$-regular (respectively, $p$-singular) cycles of $\sigma$. Now, we consider $r(\sigma)$ and $s(\sigma)$ as elements in $\Sigma_n$.

Recall that $\nu_p(n)$ denotes the largest non-negative integer such that $p^{\nu_p(n)}$ divides $n$. 

\begin{Prop}\label{regular-singular}
Let $R$ be a field of characteristic $p\ge 0$, $\sigma\in \Sigma_n$ and $\tau\in \Sigma_m$. If $S_n(c_{\sigma},R)$ and $S_m(c_{\tau},R)$ are derived equivalent, then

$(1)$ $S_n(c_{r(\sigma)},R)$ and $S_m(c_{r(\tau)},R)$ are derived equivalent, and

 $(2)$ $S_n(c_{s(\sigma)},R)$ and $S_m(c_{s(\tau)},R)$ are derived equivalent.
\end{Prop}

{\it Proof.} By Theorem \ref{main1}(1), we show the following: If $c_{\sigma}\stackrel{D}\sim c_{\tau}$, then $c_{r(\sigma)}\stackrel{D}\sim c_{r(\tau)}$ and $c_{s(\sigma)}\stackrel{D}\sim c_{s(\tau)}$.

Indeed, let $\lambda=(\lambda_1,\cdots,\lambda_k)$ be the cycle type of $\sigma$.  For $i\in [k]$ and an irreducible factor $f(x)$ of
$x^{\lambda_i}-1$, we define $q_{f(x)}:= max\{\nu_{p}(\lambda_j) \mid j\in[k], f(x)\le x^{\lambda_j}-1\}$. Then it follows from Lemma \ref{per} that $\mathcal{E}_{c_{\sigma}} = \{ f(x)^{p^{\nu_p(\lambda_i)}}\mid i\in[k], f(x)\mbox{ is irreducible and } f(x)\le x^{\lambda_i}-1\}$ and $$\mathcal{M}_{c_\sigma}=\{f(x)^{p^{q_{f(x)}}}\mid i\in [k], f(x)\mbox{ is irreducible and } f(x)\le x^{\lambda_i}-1 \}.$$ This shows that $\mathcal{M}_{c_{r(\sigma)}}=\{f(x)\mid f(x)\in \mathcal{E}_{c_\sigma}~\mbox{ is irreducible}\}$ and $\mathcal{M}_{c_{s(\sigma)}}=\{g(x)\mid g(x)\in \mathcal{M}_{c_\sigma}~\mbox{ is reducible }\}$.

Let $a$ denote the order of $r(\sigma)$. Then it follows from the definition of $r(\sigma)$ that $p\nmid a$ and therefore $m_{r(\sigma)}(x)\leq x^a-1$ is a product of distinct irreducible polynomials. Thus $\mathcal{E}_{c_{r(\sigma)}}=\mathcal{M}_{c_{r(\sigma)}}$, consisting only of some irreducible polynomials. By the definition of $s(\sigma)$, we see that $\mathcal{E}_{c_\sigma}=(\mathcal{E}_{c_{r(\sigma)}}\setminus \{x-1\})\cup \mathcal{E}_{c_{s(\sigma)}}$ and that polynomials in $\mathcal{E}_{c_{s(\sigma)}}\setminus \{x-1\}$ are reducible. Hence
$$\mathcal{E}_{c_{s(\sigma)}}=\begin{cases}\{u(x)\in \mathcal{E}_{c_\sigma}\mid u(x)~\mbox{is reducible in}~ R[x]\} & \mbox{ if } r(\sigma)=1,\\ \{u(x)\in \mathcal{E}_{c_\sigma}\mid u(x)~\mbox{is reducible in}~ R[x]\}\cup \{x-1\} & \mbox{ if } r(\sigma)\neq 1.\end{cases}$$

Suppose $c_{\sigma}\stackrel{D}\sim c_{\tau}$ . Then there is a bijection $\pi: \mathcal{M}_{c_\sigma}\to \mathcal{M}_{c_\tau}$ such that $R[x]/(h(x))\simeq R[x]/((h(x)\pi)$ as algebras and $P_{c_\sigma}(h(x))=P_{c_\tau}((h(x))\pi)$ for $h(x)\in \mathcal{M}_{c_\sigma}$. For irreducible polynomials $w(x),z(x)\in R[x]$, if $R[x]/(w(x)^a)\simeq R[x]/(z(x)^b)$ as algebras, then $a=b$ and $R[x]/(w(x)^e)\simeq R[x]/(z(x)^e)$ as algebras for all $e\le a$. Thus we may extend $\pi$ to a bijection between $\mathcal{E}_{c_\sigma}$ and $\mathcal{E}_{c_\tau}$ such that $R[x]/(h(x))\simeq R[x]/((h(x)\pi)$ as algebras for $h(x)\in \mathcal{E}_{c_\sigma}$. Since $\mathcal{E}_{c_\sigma}=(\mathcal{E}_{c_{s(\sigma)}}\setminus \{x-1\})\dot\cup \mathcal{E}_{c_{r(\sigma)}}$, the restriction of $\pi$ on $\mathcal{M}_{c_{r(\sigma)}}$ (respectively, $\mathcal{M}_{c_{s(\sigma)}}$) maps onto $\mathcal{M}_{c_{r(\tau)}}$ (respectively, $\mathcal{M}_{c_{s(\tau)}}$). For $v(x)\in \mathcal{M}_{c_{r(\sigma)}}$, there holds $P_{c_{r(\sigma)}}(v(x))=P_{c_{r(\tau)}}((v(x))\pi)=\{1\}$. Thus $c_{r(\sigma)}$ and $c_{r(\tau)}$ are $D$-equivalent.

Particularly, $r(\sigma)=1$ if and only if $r(\tau)=1$. This yields that $x-1\in \mathcal{E}_{c_{s(\sigma)}}$ if and only if $x-1\in \mathcal{E}_{c_{s(\tau)}}$. Let $a$ and $b$ be nonnegative integers such that $(x-1)^{p^a}\in \mathcal{M}_{c_\sigma}$ and $(x-1)^{p^b}\in \mathcal{M}_{c_\tau}$ are the only polynomials divisible by $x-1$, and that $((x-1)^{p^a})\pi=j(x)$ and $((x-1)^{p^b})\pi^{-1}=k(x)$. Then it follows from $R[x]/(h(x))\simeq R[x]/((h(x))\pi)$ for $h(x)\in \mathcal{M}_{c_\sigma}$ that $j(x)=(x+u)^{p^a}$ and $k(x)=(x+v)^{p^b}$ for $u,v\in R$. Since $P_{c_\sigma}(k(x))\subseteq P_{c_\sigma}((x-1)^{p^a})$ and $P_{c_\tau}(j(x))\subseteq P_{c_\tau}((x-1)^{p^b})$, we obtain
$$P_{c_\sigma}(k(x))=P_{c_\sigma}((x-1)^{p^a})=P_{c_\tau}(j(x))=P_{c_\tau}((x-1)^{p^b}).$$
Hence $a=b$. We may assume that $\pi$ maps $(x-1)^{p^a}$ in $\mathcal{M}_{c_\sigma}$ to $(x-1)^{p^b}$ in $\mathcal{M}_{c_\tau}$.
By the above calculations of $\mathcal{M}_{c_{s(\sigma)}}$ and $\mathcal{E}_{c_{s(\sigma)}}$, we get $P_{c_{s(\sigma)}}(h(x))=P_{c_\sigma}(h(x))\setminus\{1\}$ for $h(x)\in \mathcal{M}_{c_{s(\sigma)}}$ with $h(x)\neq (x-1)^{p^a}$. For $(x-1)^{p^a}$ in $\mathcal{M}_{c_\sigma}$, we have $P_{c_\sigma}((x-1)^{p^a})=P_{c_{s(\sigma)}}((x-1)^{p^a})$. Thus the matrices $c_{s(\sigma)}$ and $c_{s(\tau)}$ are $D$-equivalent. $\square$

\medskip
For a counterexample to the converse of Proposition \ref{regular-singular}, we refer to Example \ref{ex4.5} in the next section.

\section{Examples and questions\label{sect5}}
In this section we display a few examples to illustrate our results in the previous sections.

The following example shows that the centralizer matrix algebras of non-conjugate matrices may be Morita equivalent.

\begin{Bsp}{\rm
Let $R$ be a field and $J_n(\lambda)$ the $n\times n$ Jordan matrix with the eigenvalue $\lambda\in R$. We take $c=J_3(1)\oplus J_4(1)\oplus J_3(0)\oplus J_2(0)$ and $d=J_3(0)\oplus J_4(0)\oplus J_3(1)\oplus J_2(1)$. In general, we have $m_{c\oplus d}(x)=[m_c(x),m_d(x)]$, where $[f(x),g(x)]$ stands for the least common multiple of $f(x)$ and $g(x)$ in $R[x]$. Then $m_c(x)=x^3(x-1)^4$, $\mathcal{E}_c=\{x^2, x^3, (x-1)^3, (x-1)^4\}$, $\mathcal{M}_c=\{x^3,(x-1)^4\}$, $P_c(x^3)=\{2,3\},P_c((x-1)^4)=\{3,4\}$, and $m_d(x)=x^4(x-1)^3$, $\mathcal{E}_d=\{x^3,x^4,(x-1)^2,(x-1)^3\}$, $\mathcal{M}_d=\{x^4, (x-1)^3\}$,$P_d(x^4)=\{3,4\}$, $P_d((x-1)^3)=\{2,3\}$. Let $\pi: \mathcal{M}_c\to \mathcal{M}_d$ be the map: $x^3\mapsto (x-1)^3, (x-1)^4\mapsto x^4$. Then it follows from Theorem \ref{main1}(1) that $S_{12}(c,R)$ and $S_{12}(d,R)$ are Morita equivalent, while $c$ and $d$ are not conjugate since they have different minimal polynomials.}
\end{Bsp}

The next example shows that the existence of a Morita equivalence between principal centralizer matrix algebras depends on the ground field.

\begin{Bsp} {\rm Let $\sigma:=(1~2~3~4~5)(6~7~8\cdots~17~18), \tau:=(1~2~3~4~5~6~7)(8~9\cdots ~17~18)\in \Sigma_{18}$. The minimal polynomials of $c_\sigma$ and $c_\tau$ over $\mathbb{Q}$ are $(x-1)(x^4+x^3+x^2+x+1)(x^{12}+x^{11}+\cdots +x+1)$ and $(x-1)(x^{10}+x^9+\cdots+x+1)(x^6+x^5+\cdots +x+1)$, respectively. In this case, $\mathcal{M}_{c_{\sigma}}=\{x-1, x^4+x^3+x^2+x+1, x^{12}+x^{11}+\cdots +x+1\}$ and $\mathcal{M}_{c_{\tau}}=\{x-1, x^{10}+x^9+\cdots+x+1, x^6+x^5+\cdots +x+1\}$. By Theorem $\ref{main1}(1)$, $S_{18}(c_\sigma, \mathbb{Q})$ and $S_{18}(c_\tau, \mathbb{Q})$ are not Morita equivalent, while $S_{18}(c_\sigma, \mathbb{C})$ and $S_{18}(c_\tau, \mathbb{C})$ are Morita equivalent (see also \cite[Theorem 1.2(2)]{xz2}).}
\end{Bsp}

Now, we show that even for centralizer matrix algebras, almost $\nu$-stable derived equivalences may not always arise from Morita equivalences.

\begin{Bsp}{\rm
Let $a=J_5(0)\oplus J_4(0)\oplus J_2(0)$ and $b=J_5(0)\oplus J_3(0)\oplus J_1(0)$. Then the centralizer algebras $S_{11}(a,R)$ and $S_9(b,R)$ are not Morita equivalent, but they are almost $\nu$-stable derived equivalent by Theorem \ref{main1}(3).}
\end{Bsp}

We point out that even in the class of principal centralizer matrix algebras, derived equivalences do not have to preserve representation-finiteness, while almost $\nu$-stable derived equivalences always preserve representation-finiteness for arbitrary algebras.

\begin{Bsp}{\rm
Let $R$ be an algebraically closed field, $c:=J_5(0)\oplus J_4(0)\oplus J_1(0)\in M_{10}(R)$ and $d:=J_5(0)\oplus J_2(0)\oplus J_1(0)\in M_8(R)$. Then $S_{10}(c,R)$ and $S_8(d,R)$ are derived equivalent by Theorem \ref{main1}(2), while $S_{10}(c,R)$ is representation-finite, but $S_8(d,R)$ is not by Lemma \ref{rep-f}.
}\end{Bsp}

\medskip
Generally, the converse of Proposition \ref{regular-singular} may be false, as we can see by the following example.

\begin{Bsp}\label{ex4.5} {\rm Let $R$ be an algebraically closed field of characteristic $5$. We take $\sigma\in \Sigma_{19}$ with the cycle type $(15,4)$, and $\tau\in \Sigma_{20}$ with the cycle type $(15,3,2)$. In this case, $r(\sigma)$ is a permutation of the cycle type $(4,1^{15})$ and $s(\sigma)$ is a permutation of cycle type $(15, 1^4)$, while $r(\tau)$ has the cycle type $(3,2, 1^{15})$ and $s(\tau)$ has the cycle type $(15,1^5)$. Clearly,  $S_{19}(c_{s(\sigma)},R)$ and $S_{20}(c_{s(\tau)},R)$ are derived equivalent by Corollary \ref{one-more}. Since
$\mathcal{M}_{c_{r(\sigma)}}= \{x-1, x+1, x-\eta, x+\eta\}$ and $\mathcal{M}_{c_{r(\tau)}}= \{x-1, x+1, x+\epsilon, x-\epsilon^2\}$, where $\eta$ and $\epsilon$ are $4$-th and $3$-th primitive roots of unity, respectively, it follows from Theorem \ref{main1}(2) that
$S_{19}(c_{r(\sigma)},R)$ and $S_{20}(c_{r(\tau)},R)$ are derived equivalent.

By Lemma \ref{per}, $\mathcal{M}_{c_\sigma}=\{(x-1)^5,(x-\epsilon)^5,(x-\epsilon^2)^5,x+1,x-\eta,x+\eta\}$ and $\mathcal{M}_{c_\tau}=\{(x-1)^5,(x-\epsilon)^5,(x-\epsilon^2)^5,x+1\}$. Clearly, $|\mathcal{M}_{c_\sigma}|=6\neq 4=|\mathcal{M}_{c_\tau}|$. Hence no bijections between $\mathcal{M}_{c_\sigma}$ and $\mathcal{M}_{c_\sigma}$ exist, and therefore $S_{19}(c_{\sigma},R)$ and $S_{20}(c_{\tau},R)$ cannot be derived equivalent by Theorem \ref{main1}.

This shows that derived equivalences of the centralizer matrix algebras of $p$-regular parts and $p$-singular parts of permutations do not have to guarantee a derived equivalence of the ones of the permutations themselves.
}\end{Bsp}

Having described derived equivalences of principal centralizer matrix algebras, we propose the following questions for further study. In the following, $R$ stands for a field.

{\bf Question 1}. Under which necessary and sufficient conditions on matrices $c\in M_n(R)$ and $d\in M_m(R)$ are $S_n(c,d)$ and $S_m(d,R)$ stably equivalent?

{\bf Question 2}. For which permutations $\sigma\in \Sigma_n$ and $\tau\in \Sigma_m$ do derived equivalences between $S_n(c_{s(\sigma)},R)$ and $S_m(c_{s(\tau)},R)$, and between $S_n(c_{r(\sigma)},R)$ and $S_m(c_{r(\tau)},R)$ ensure a derived equivalence between $S_n(c_{\sigma},R)$ and $S_m(c_{\tau},R)$ ?

Related to general consideration of the centralizers of matrices, we mention the following.

{\bf Question 3}. Describe structural and homological properties of $S_n(C,R)$ for $|C|\ge 2$.

\medskip
\textbf{Acknowledgements.} The research work was supported partially by the National Natural Science Foundation of China (Grants 12031014 and 12226314). The authors thank Dr. Jinbi Zhang for discussions on the primary version of the manuscript.

{\footnotesize
\smallskip
Xiaogang Li,
School of Mathematical Sciences, Capital Normal University, 100048
Beijing, P. R. China;

{\tt Email: 2200501002@cnu.edu.cn}

\smallskip
Changchang Xi,

School of Mathematical Sciences, Capital Normal University, 100048
Beijing; \&
School of Mathematics and Statistics, Shaanxi Normal University, 710119 Xi'an, P. R. China

{\tt Email: xicc@cnu.edu.cn}
}


{\footnotesize
\begin{thebibliography}{99}
\bibitem{Al}{{\sc J. L. Alperin}, Weights for finite groups, \emph{Proc. Symposia Pure Math.} \textbf{47} (1987) 369-379.}

\bibitem{Al2}{{\sc J. L. Alperin}, A Lie approach to finite groups, In: \emph{Lecture Notes in Mathematics} \textbf{1456}, Groups-Canberra Sringer-Verlag New York, (1989) 1-8.}

\bibitem{As}{{\sc H. Asashiba}, On a lift of an individual stable equivalence to a standard derived equivalence for representation-finite self-injective algebras, \emph{Algebr. Represent. Theory} \textbf{6} (2003) 427-447.}

\bibitem{ARS}{{\sc M. Auslander}, {\sc I. Reiten} and {\sc S. O. Smal\o}, {\it Representation of artin algebras}, Cambridge Studies in Advanced Mathematics \textbf{36}, Cambridge University Press, Cambridge, 1995.}

\bibitem{bf}{{\sc L. Brickman} and {\sc P. A. Fillmore}, The invariant subspace lattice of a linear
transformation, \emph{Canad. J. Math.} \textbf{19} (1967) 810-822.}

\bibitem{MB}{{\sc M. Brou\'e}, Equivalences of blocks of group algebras, In: {\it Finite-dimensional algebras and related topics} (Ottawa, ON, 1992), 1-26, NATO Adv. Sci. Inst. Ser. C: Math. Phys. Sci. 424, Kluwer Acad. Publ., Dordrecht, 1994.}

\bibitem{ca}{{\sc M. Cabanes}, Brauer morphism between modular Hecke algebras, \emph{J. Algebra} \textbf{115} (1988) 1-31.}

\bibitem{CM}{{\sc A. Chan} and {\sc R. Marczinzik}, On representation-finite gendo-symmetric biserial algebras, \emph{Algebr. Represent. Theory} \textbf{22} (1) (2019) 141-176.}

\bibitem{cdfk}{{\sc H. F. da Cruz, G. Dolinar, R. Fernandes} and {\sc B. Kuzma}, Maximal doubly stochastic matrix centralizers, \emph{Linear Algebra Appl.} \textbf{532} (2017) 387-398.}

\bibitem{dr}{{\sc V. Dlab} and {\sc C. M. Ringel}, The module theoretical approach to quasi-hereditary algebras,
London Math. Soc. Lect. Note Ser. \textbf{168} (1992) 200-224.}

\bibitem{YV}{{\sc Y. Drozd} and {\sc V. Mazorchuk}, Representation type of $\; {}_{\lambda}^{\infty}{\mathcal{H}^1_{\mu}}$, \emph{Quart. J. Math.} \textbf{57} (2006) 319-338.}

\bibitem{hr}{{\sc L. Hille} and {\sc G. R\"ohrle}, A classification of parabolic subgroups of classical groups with a finite number
of orbits on the unipotent radical, \emph{Transform. Groups} \textbf{4} (1) (1999) 35-52.}

\bibitem{hx1}{{\sc W. Hu} and {\sc C. C. Xi}, Derived equivalences and stable equivalences of Morita type I, \emph{Nagoya Math. J.} \textbf{200} (2010) 107-152.}

\bibitem{hx3}{{\sc W. Hu} and {\sc C. C. Xi}, Derived equivalences and stable equivalences of Morita type II, \emph{Rev. Mat. Iberoam.} \textbf{34} (1) (2018) 59-110.}
\bibitem{hx2}{{\sc W. Hu} and {\sc C. C. Xi}, $\mathcal{D}$-split sequences and derived equivalences, \emph{Adv. Math.} \textbf{227} (1) (2011) 292-318.}

\bibitem{JL}{{\sc C. U. Jensen} and {\sc H. Lenzing}, Homological dimension and representation type of algebras under base field extension, \emph{Manuscr. Math.} \textbf{39} (1982) 1-13.}

\bibitem{kaplansky}{{\sc I. Kaplansky}, {\it Linear algebra and geometry}, a second course, Chelsea Publishing Company, New York, 1974.}

\bibitem{Liu1}{{\sc Y. M. Liu}, Summands of stable squivalences of Morita type,
\emph{Comm. Algebra} \textbf{36} (10) (2008) 3778-3782.}

\bibitem{Liu2}{{\sc Y. M. Liu} and {\sc C. C. Xi}, Constructions of stable equivalences of Morita type
for finite dimensional algebras II, \emph{Math. Z.} \textbf{251} (2005) 21-39.}

\bibitem{LX3}{{\sc Y. M. Liu} and {\sc C. C. Xi}, Constructions of stable equivalences of Morita type
for finite dimensional algebras III, \emph{J. Lond. Math. Soc.} \textbf{76} (2)(2007) 567-585.}

\bibitem{MV1}{{\sc R. Mart\'{i}nez-Villa}, Algebra stably equivalent to $l$-hereditary, In: {\it Representation Theory,} II (Proc. Second Internet. Conf., Carleton Univ., Ottawa, Ont., 1979), pp. 396-431, Lecture Notes in Math. \textbf{832}, Springer, Berlin, 1980.}

\bibitem{MV2}{{\sc R. Mart\'{i}nez-Villa}, Properties that are left invariant under stable equivalence, \emph{Comm. Algebra} \textbf{18} (12) (1990) 4141-4169.}

\bibitem{BJ}{{\sc B. J. M\"{u}ller}, The classification of algebras by dominant dimension, \emph{Canad. J. Math.} \textbf{20} (1968) 398-409.}

\bibitem{DP}{{\sc D. I. Panyushev}, Two results on centralisers of nilpotent elements, \emph{J. Pure Appl. Algebra} \textbf{212} (2008) 774-779.}

\bibitem{AP}{{\sc A. Premet}, Nilpotent commuting varieties of reductive Lie algebras, \emph{Invent. Math.} \textbf{154} (2003) 653-683.}

\bibitem{IR}{{\sc I. Reiten}, Stable equivalence of self-injective algebras, \emph{J. Algebra} \textbf{40} (1) (1976) 60-74.}

\bibitem{Rickard1}{{\sc J. Rickard}, Morita theory for derived categories, \emph{J. Lond. Math. Soc.} \textbf{39} (2) (1989) 436-456.}

\bibitem{JR2}{{\sc J. Rickard}, Derived equivalences as derived functors, \emph{J. Lond. Math. Soc.} \textbf{43} (1) (1991) 37-48.}

\bibitem{rotman}{{\sc J. J. Rotman}, \emph{An introduction to homological algebra}, Pure and Applied Mathematics \textbf{85}, Academic Press,  New York-London, 1979.}

\bibitem{rouq}{{\sc R. Rouquier}, Derived equivalences and finite dimensional algebras. In: \emph{International Congress of
Mathematicians}, Vol. II, 191-221. Eur. Math. Soc., Zurich, 2006.}

\bibitem{HW}{{\sc H. Weyl}, {\it The classical groups. Their invariants and representations}, Princeton University Press, Princeton, N.J., 1939.}

\bibitem{We}{{\sc J. H. M. Wedderburn}, On hypercomplex numbers, \emph{Proc. Lond. Math. Soc.} \textbf{6} (1908) 77-118.}

\bibitem{x3}{{\sc C. C. Xi}, Derived equivalences of algebras, \emph{Bull. Lond. Math. Soc.} \textbf{50} (6) (2018) 945-985.}

\bibitem{xz1}{{\sc C. C. Xi} and {\sc J. B. Zhang}, Structure of centralizer matrix algebras, \emph{Linear Algebra Appl.} \textbf{622} (2021) 215-249.}

\bibitem{xz2}{{\sc C. C. Xi} and {\sc J. B. Zhang}, Centralizer matrix algebras and symmetric polynomial of partitions, \emph{J. Algebra} \textbf{609} (2022) 688-717.}

\bibitem{xz3}{{\sc C. C. Xi} and {\sc J. B. Zhang}, New invariants of stable equivalences of algebras, Preprint, 1-20. arXiv:2207.10848.}

\bibitem{Z}{{\sc A. Zimmermann}, {\it Representation theory. A homological algebra point of view}, Algebra and Application \textbf{19}, Springer Cham, 2014.}

\end{thebibliography}
}
\end{document}